# Scalability and Performance Trade-offs Between Q-Learning and Parallel Intensified Simulated Annealing in Flexible Job Shop Scheduling

Robert L Burdett (r.burdett@qut.edu.au) and Paul Corry (p.corry@qut.edu.au)
School of Mathematical Sciences, Queensland University of Technology, Brisbane, Qld, Australia

**ABSTRACT**

This article presents a controlled and reproducible comparative study of tabular Q-learning (QL) and Parallelised Simulated Annealing with Intensification (PISA) for the Flexible Job Shop Scheduling Problem (FJSP). Although reinforcement learning and metaheuristics have both been widely studied in scheduling, direct comparisons are often difficult due to inconsistent modelling assumptions, representations, and computational settings. To address this gap, the study develops a unified open implementation framework that supports both deterministic and stochastic FJSP environments under consistent evaluation budgets and solution representations. Within the Q-learning framework, methodological refinements are introduced through Pareto-based action selection and a mean–variance Pareto extension designed to balance solution quality and robustness under uncertainty.

The study conducts a controlled comparison against a strong intensified metaheuristic baseline using the Dauzère-Pérès and Paulli benchmark instances. Attention is given to practical factors that substantially affect performance but are frequently under-reported, including simulation budget, state-space growth, memory requirements, and the absence of parallelism in standard tabular Q-learning. The analysis isolates the conditions under which tabular Q-learning remains competitive and identifies where structural limitations become dominant.

Results show that PISA generally provides superior computational efficiency, scalability, and consistency in best-case performance across the benchmark set. Tabular Q-learning can achieve competitive results on selected small and medium-scale instances, particularly when sufficient training effort is available, but its performance deteriorates as state-space complexity and sampling requirements increase. The findings indicate that memory growth, sample efficiency, and sequential learning behaviour constitute the primary practical constraints limiting the applicability of tabular Q-learning for larger FJSP instances. Overall, the study clarifies the trade-offs between learning-based and metaheuristic approaches and provides practical guidance regarding when tabular Q-learning is likely to be beneficial in scheduling environments and when intensified metaheuristics remain the more effective choice.

## 1. Introduction

The Flexible Job Shop Scheduling Problem (FJSP) remains one of the most challenging combinatorial optimisation problems due to its coupled assignment and sequencing decisions. In the FJSP, a set of jobs $J$, each consisting of an ordered sequence of operations $O_j = \{o_{j,k}\}$ for $j \in J$, must be scheduled on a set of machines $M$. Unlike the classical job shop problem (JSP), each operation may be processed on multiple alternative machines, with processing times $p_{j,k}$ that can vary substantially across machines. This added flexibility significantly increases the complexity of decision-making (Gomez et al. 2023).

The objective of the FJSP problem is to construct feasible schedules that minimise makespan or related performance measures subject to standard operational constraints, including machine capacity, job precedence relations, and operation-specific requirements. The problem is well known to be NP-hard, with two principal sources of difficulty: (i) exponential growth of the joint decision space as the number of jobs and machines increases, and (ii) uncertainty introduced by stochastic or variable processing times, which weakens the assumptions underlying deterministic scheduling models. Consequently, exact optimisation methods are not computationally viable for realistically sized instances, motivating the use of heuristic and learning-based approaches.

Despite extensive recent activity in reinforcement learning (RL) for scheduling (Koszykowski and Orzeszko (2026), Nair and Bhanu (2022)), evaluation practices remain inconsistent and evidence supporting comparative effectiveness across operational settings is still evolving (Dulac-Arnold et al. (2019), Srinivasan (2023), Kommey et al. (2024)). While these methods show considerable promise, prior studies have reported practical limitations and instances of inconsistent performance (Shao et al. 2025). Many Q-learning-based FJSP studies rely on restricted episode budgets and do not report key indicators such as state-space growth, memory consumption, or convergence behaviour. This limits the ability to assess scalability and hinders reliable cross-study comparison. Moreover, claims regarding the superiority of RL in stochastic or dynamic scheduling settings are not always supported by evaluations against strong, carefully implemented metaheuristic baselines. These observations motivate the following research questions:

RQ1: Under matched computational budgets, where does tabular Q-learning sit relative to a strong meta-heuristic baseline?
RQ2: How do state-space growth and episode budget shape Q-learning performance and memory footprint?
RQ3: For the FJSP does effective application of tabular Q-learning depend on problem-specific state representations, design choices, and informed action-selection mechanisms
RQ4: Does Pareto-based action selection improve on ε-greedy for the FJSP, and by how much?
RQ5: How do these conclusions change under stochastic processing times?

These hypotheses are tested empirically in this study. To this end, we implement two solution approaches: (i) a tabular Q-learning algorithm and (ii) a metaheuristic based on Simulated Annealing with Intensification (ISA). Both are implemented using comparable levels of design and coding effort to ensure a balanced comparison. The study considers both deterministic and stochastic processing times, as this distinction materially affects algorithm behaviour and performance interpretation. In the stochastic setting, the objective is defined as minimisation of expected makespan, thereby focusing on average-case performance without incorporating explicit robustness or risk measures to mitigate worst-case outcomes. Those considerations are left for future work.

We focus on tabular Q-learning as the representative RL method to isolate the contribution of state representation, memory growth and exploration strategy from neural function approximation. Recent DRL work conflates these design choices; the controlled setting adopted here allows their effects to be measured independently.

For the metaheuristic baseline, we employ Simulated Annealing with Intensification (ISA), a strong and well-established hybrid variant of SA. SA remains widely used in combinatorial optimisation and continues to serve as a competitive baseline when carefully implemented (Burdett and Kozan (2010, 2014, 2018), Burdett et al. (2019), Araujo et al (2025)). Carefully designed SA implementations have frequently been shown to match or outperform more complex hybridised approaches, which motivates its use here as a robust benchmark for evaluating Q-learning performance. We restrict the comparison to well-established optimisation frameworks with documented empirical performance, rather than attempting an exhaustive survey across a broader class of heuristic methods. Newer metaphor-based metaheuristics for instance have been excluded as their development has often prioritised conceptual novelty over consistent empirical validation and reproducible effectiveness (Rolim and Nagano, 2025).

The contributions of this paper are fourfold:

1. To conduct a controlled and reproducible comparison between tabular Q-learning and ISA for the FJSP under similar evaluation budgets.

2. To provide a consistent experimental assessment against a strong metaheuristic baseline.
3. To analyse implementation-level factors affecting Q-learning performance, including memory growth and sensitivity to simulation budget.
4. To identify problem regimes in which tabular Q-learning is competitive and those in which structural limitations dominate.

Implementation details are documented explicitly to support reproducibility and to highlight design choices that materially influence algorithmic performance. The full implementation, including SA and Q-learning frameworks, is publicly available at: https://github.com/rburdett76/QL_PISA.

The remainder of the paper is organised as follows. Section 2 reviews relevant literature. Section 3 presents the proposed SA and Q-learning methodologies. Section 4 reports computational results and analysis. Section 5 concludes the study and outlines directions for future research.

## 2. Literature Review

This section reviews existing literature relevant to the problem domain, highlighting key approaches, findings, and gaps.

### 2.1. Reinforcement Learning.

Numerous RL approaches have been applied to job scheduling problems. Table 1 summarises a selection of recent representative studies, although it is not intended to be exhaustive. A substantial body of work on Q-learning can be found in the earlier literature (circa 2005–2017) (Gui et al. 2023), but comparatively fewer studies have appeared in more recent years. For a more comprehensive overview of RL applications in this domain, readers are referred to the review articles in Zhang and Zhu (2025) and Rihane et al. (2025).

**Table 1.** Recent studies on RL for JSP and FJSP

| Article (Year) | Prob. | Stoch. Proc. Time | Dyn. Jobs | Obj. | Approach | Comparative. Approach |
|---|---|---|---|---|---|---|
| Fonseca-Reyna et al. (2017) | FSSP | × | × | Cmax | QL | Cited benchmarks |
| Lang et al. (2020) | FJSP | × | × | Cmax | DQN | GRASP |
| Luo (2020) | FJSP | × | ✓ | Tardiness | DQN | DR, QL |
| Han and Yang (2021) | FJSP | × | × | Cmax | DRL+RNN | DR |
| Luo et al (2021) | FJSP | × | ✓ | Tardiness | DQN | DR, QL |
| Palacio et al. (2022) | FJSP | × | × | Cmax | Two-stage QL | - |
| Belmamoune et al. (2022) | JSP | × | × | Cmax | QL | GA, DRL, DR |
| Zhou et al. (2022) | JSP | × | ✓ | Composite | RL | GA, DR |
| Gui et al. (2023) | FJSP | × | ✓ | Tardiness | DQN | DR |
| Zhang et al. (2023) | FJSP | × | × | Cmax | MA_DRL+DQN | DR |
| Smit et al. (2024) | FJSP | ✓ | × | E[Cmax] or VaR | DRL + SPM-DAN | CP, DR |
| Tiacci et al. (2024) | FJSP | × | ✓ | Tardiness | DRL + DNN | DR |
| Wang et al. (2024) | FJSP | × | × | Cmax | DRL + DAN | GA, DR |
| Infantes et al. (2024) | FJSP | ✓ | ✓ | E[Cmax] | DRL + GNN | DR, OR-Tools, CP |
| Shao et al. (2025) | FJSP | × | × | Cmax | DRL + GNN | DR, OR-Tools |
| Huang et al. (2025) | FJSP | × | × | Cmax | DRL + GNN | GA |
| Liu et al. (2025) | FJSP | × | × | Cmax | DRL + GNN | GA, DR, OR-Tools |

Key. ABC (Artificial Bee Colony), AC (Actor-Critic), Cmax (Makespan), CP (Constraint Programming), DAN (Dual attention network), DR (Dispatching Rules), DQN (Deep Q-Networks), DNN (Deep Neural Network), DRL (Deep Reinforcement Learning), FSSP (flow shop scheduling problem), GA (Genetic Algorithm), GNN (Graph Neural Network), MA_DRL (Multi-Agent DRL), NCO (Neural Combinatorial Optimization), PPO (Proximal Policy Optimization), SPM (Scenario Processing Module), RNN (Recurrent Neural Network), VaR (Value at Risk)

In our review, we observed that comparatively few studies include direct benchmarks against established metaheuristics such as Simulated Annealing. We also note that a substantial portion of recent research on RL applied to the job shop and flexible job shop scheduling problems originates from China, reflecting strong regional interest and investment in smart manufacturing and intelligent optimisation.

Much of the existing literature focuses on specialised or modified problem variants, with fewer studies systematically evaluating performance on standard benchmark formulations. In addition, reporting practices vary considerably across studies; for example, conclusions are often presented in broad terms such as "outperforms method X" or "is feasible," while details such as computational effort, runtime, and the practical significance of observed improvements are not always fully reported. These differences in emphasis and reporting reduce the ease of reproducibility and cross-study comparison, and they suggest an opportunity for future work to adopt more consistent evaluation protocols and to provide fuller disclosure of implementation and performance details.

Srinivasan (2023), Kommey et al. (2024) and Rihane et al. (2025) provide a comprehensive discussion of key limitations associated with RL methods. As highlighted in these studies, RL algorithms typically learn through interaction with simulated or "virtual" environments, often requiring many interactions, which makes sample efficiency a central challenge.

A commonly reported issue is limited scalability, where performance can deteriorate as problem size and complexity increase. Another fundamental difficulty lies in the exploration–exploitation trade-off: agents must balance exploring the search space to identify promising strategies while exploiting current knowledge to improve performance.

RL methods are also frequently affected by sparse, delayed, or noisy reward signals, which complicate the attribution of actions to long-term outcomes. In addition, policies learned in specific training environments may exhibit limited generalisation, with performance degrading when transferred to new or modified problem instances.

A comprehensive analysis of the strengths and limitations of the articles presented in Tables 1 is provided below:

Fonseca-Reyna et al. (2017) considers flow shop scheduling, where the QL approach defines states as partial job sequences and selects the next job from the set of currently available options. The proposed method is evaluated using Taillard benchmark instances. The reported results are promising for small- to medium-sized instances, with the QL approach generating competitive solutions in these settings. However, performance appears to deteriorate on larger instances. In addition, solution times are not reported, which limits the ability to assess computational efficiency relative to state-of-the-art metaheuristic approaches.

Lang et al. (2020) investigates the integration of Deep Q-Networks (DQN) with discrete event simulation (DES) for a variant of the FJSP. Their numerical experiments suggest that high-quality solutions can be obtained with relatively limited training in the tested settings. However, several aspects of this line of work indicate potential limitations. The study is restricted to deterministic processing times and relatively small benchmark instances (5–8 machines and 5–20 jobs), which limits conclusions about scalability. In addition, reported computational times are higher than those observed for a greedy randomized adaptive search procedure (GRASP) baseline. The method also appears sensitive to parameter tuning, which may complicate practical implementation. Finally, the authors note that DQN approaches may be less well suited to sequencing problems due to the requirement of a fixed action space, a constraint that can limit flexibility in scheduling contexts.

Luo (2020) addresses a dynamic FJSP, with a particular focus on the insertion of new jobs. The study adopts a DRL approach to minimise total tardiness. At each rescheduling point, the DRL agent selects an appropriate dispatching rule based on the current state of the system. Building on this work, Luo et al. (2021) extends the framework to a multi-objective setting, considering criteria such as total weighted tardiness and machine utilisation. The authors propose a hierarchical architecture consisting of two Deep Q-Networks (DQNs): a higher-level DQN that determines the current optimisation objective, and a lower-level DQN that acts as a decision-making module. The lower-level DQN receives both state features and the objective specified by the higher-level network and selects an appropriate dispatching rule accordingly.

Palacio et al. (2022) implemented a RL algorithm to tackle a real-world flexible job scheduling problem in a manufacturing setting. Their approach is notable for its two-phase structure, both leveraging QL. The first phase focuses on learning the most suitable machine assignments for operations, while the second determines the optimal execution order of those operations. While the idea shows promise, the lack of supporting numerical data indicates an opportunity for future empirical study.

Belmamoune et al. (2022) applies a QL algorithm to solve the JSP. In their implementation, an agent is trained to select among four dispatching rules, namely first-in-first-out (FIFO), last-in-first-out (LIFO), shortest processing time (SPT), and longest processing time (LPT). The proposed approach is evaluated against a GA and a DRL method, demonstrating competitive performance. However, while the reported results are encouraging, they do not establish QL as consistently superior to alternative approaches. Further investigation is therefore required to assess its robustness and general effectiveness across a broader range of problem instances and settings.

Zhou et al. (2022) implemented a RL algorithm to address a JSP involving dynamically arriving jobs in a manufacturing context. While the study reports performance improvements over traditional heuristics, it provides limited methodological detail, making replication and in-depth evaluation challenging. Nonetheless, the work highlights the potential of RL in dynamic scheduling scenarios, even if further clarity and transparency would enhance its practical utility.

Tiacci and Rossi (2024) addresses the challenges of applying RL to the dynamic FJSP in the absence of a dedicated simulation environment. To this end, the authors develop an agent-based discrete-event simulator and integrate it with a DRL framework. While the study is methodologically innovative, its primary focus is on demonstrating the feasibility of the proposed integration rather than providing extensive performance validation. The numerical experiments include comparisons with traditional dispatching rules; however, the scope of benchmarking and comparative analysis could be further expanded to strengthen the empirical evaluation.

Han and Yang (2021), Zhang et al. (2023), Huang et al. (2025), Shao et al. (2025), Wang et al. (2025), Infantes et al. (2024), and Liu et al. (2025) employ related approaches that integrate actor–critic frameworks, graph neural networks, and deep reinforcement learning (DRL). While these methods are technically sophisticated, their complexity may present challenges for practical implementation, underscoring the need for more accessible and easily deployable formulations.

In addition, these approaches diverge substantially from conventional methods typically used for solving the FJSP. Smit et al. (2024) and Infantes et al. (2024) incorporate stochastic processing times, while Zhang et al. (2023) considers a multi-agent setting. However, the comparative evaluations in these studies are generally limited to constraint programming, OR-Tools, or dispatching rules, and do not include dedicated state-of-the-art metaheuristic baselines, which are commonly used as reference methods in the literature. Furthermore, Infantes et al. (2024) evaluates expected makespan using a relatively small set of fifty scenarios, which may limit the robustness of the reported conclusions.

### *2.2. Approaches Other Than RL*

Table 2 presents a complementary overview of approaches other than RL. A comprehensive analysis of the strengths and limitations of the articles presented in Tables 2 is provided below:

Chen et al. (2020) and Long et al. (2021) introduce self-learning metaheuristics for the basic FJSP, incorporating RL to adapt parameter control during the search process. In these studies, GA variants are adopted as benchmark methods, while comparisons with other established state-of-the-art metaheuristics are limited. In contrast, [31] provides a broader experimental assessment by including eight additional techniques, thereby offering a more comprehensive comparative perspective.

Caldeira and Gnanavelbabu (2020) proposes a simulation–optimisation approach for the FJSP with stochastic processing times. To estimate expected makespan, the authors employ Monte Carlo simulation. Their primary solution method is JAYA, a population-based metaheuristic whose distinguishing mechanism involves identifying the best and worst solutions in the current population and updating individuals according to their relative proximity to these extremes. One methodological consideration in the study is the use of a fixed number of simulations (e.g., 1000) for each expected makespan evaluation. The suitability of this choice depends on the underlying variability of processing times and may therefore be either excessive or insufficient under different stochastic conditions.

**Table 2.** Other approaches for FJSP

| Article (Year) | Prob. | Stoch. Proc. Time | Dyn. Jobs | Obj. | Approach | Comp. Approach |
|---|---|---|---|---|---|---|
| Caldeira and Gananavelbabu (2020) | FJSP | ✓ | × | E[Cmax] | MCS + JAYA | GA, ABC |
| Chen et al. (2020) | FJSP | × | × | Cmax | SeL-GA | SARSA, GA |
| Long et al (2021) | FJSP | × | × | Cmax | SeL-ABC | GWO, HA, GA |
| Gomez (2021) | FJSP | ✓ | × | SL[Cmax] | MCS + TS | - |
| Lin et al. (2022 | FJSP | ✓ | × | E[Cmax] | MCS | DR, GA, PSO |
| Tamssaouet et al. (2022) | FJSP | × | × | Two | AMOSA | Benchmarks |
| Gomez et al. (2023) | FJSP | ✓ | × | SL[Cmax] | MCS + TS | - |
| Li et al. (2023) | FJSP | × | × | Cmax & TEC | MA | 6 MO MH |
| Fu et al. (2025) | DFJSP | ✓ | × | E[Cmax], E[Tard] | MCS + MOEA | NSGA2, MA, SA |
| Araujo et al, (2025) | FJSP | × | × | Cmax | SA, TS, ILS, GRASP | - |
| Guan et al. (2025) | FJSBSP | × | × | Cmax, Energy | NSGA | - |

Key. AMOSA (Archived MO SA), BSP (Batch Processing), DJSP (Distributed FJSP), EA (Evolutionary Alg.), E[] (Expected value), GWO (Grey Wolf Optimizer), JAYA (Meta-heuristic), MA (Memetic Alg.), MCS (Monte Carlo Simulation), MH (Meta-heuristic), MO (Multi-objective), MOEA (MO EA), NSGA (Non-dominated GA), PSO (Particle Search), SARSA (On-policy RL), SL (Service Level), SeL (Self-learning), Tard (Total Tardiness) TES (Total Energy Consumption), TS (Tabu Search)

Gomez et al. (2021, 2023) examine the FJSP with stochastic processing times, using makespan service level as the primary performance objective. Their approach combines Tabu Search (TS) with Monte Carlo sampling and is applied to standard benchmark instances. These studies provide extensive computational experiments that highlight the importance of incorporating service-level considerations into scheduling decisions. However, the absence of additional baseline methods for comparison, together with limited implementation detail, constrains the transparency and reproducibility of the reported findings.

Lin et al. (2022) applies a Grey Wolf Optimization (GWO) meta-heuristic to the FJSP with stochastic processing times. Notably, the study incorporates a RL component to dynamically adjust algorithm parameters with the aim of improving search performance. It also introduces an Optimal Computing Budget Allocation (OCBA) mechanism, which estimates key statistical measures to manage the number of Monte Carlo simulation scenarios required for evaluation. The numerical experiments are conducted within a semiconductor manufacturing context. However, the overall experimental scope is relatively narrow, which may limit the robustness and generalisability of the performance assessment.

Li et al. (2023) proposes a Learning-Based Reference Vector Memetic Algorithm (LRVMA) for a multi-objective FJSP, with the objectives of minimising total energy consumption and makespan. Processing times are represented using type-2 fuzzy uncertainty; however, this uncertainty is ultimately transformed into a single expected value for each task through a credibility-based defuzzification procedure, resulting in a deterministic evaluation framework rather than a fully stochastic one.

The algorithm incorporates standard crossover operators, including precedence operation crossover and uniform crossover, together with local search procedures. Its principal contribution lies in the use of a QL mechanism to dynamically adjust neighbourhood size, thereby influencing solution diversity during the search process. Although the study benchmarks LRVMA against six additional metaheuristics, the presentation of comparative methods and implementation settings is limited. Algorithmic descriptions, supporting citations, parameter configurations, tuning procedures, termination criteria, and perturbation mechanisms are not reported in detail. As a result, independent replication and direct methodological comparison may be challenging based solely on the information provided.

Fu et al (2025) considers a distributed FJSP with stochastic processing times and multiple objectives. The study applies an evolutionary algorithm and benchmarks its performance against SA,

NSGA-II, and a memetic algorithm (MA). While the comparative scope is broad, the rationale for selecting this set of benchmark methods is not discussed in detail, and implementation specifics are only briefly reported. In addition, the Monte Carlo evaluation is based on 50 samples, which may limit the statistical robustness of the estimated performance measures.

Araujo et al (2025) investigated the application of local search and four trajectory meta-heuristics (i.e., Iterated Local Search (ILS), Greedy Randomized Adaptive Search (GRASP), Tabu Search (TS), and Simulated Annealing (SA)) to the FJSP with non-serial operation precedence constraints and position-based learning effects. Extensive computational experiments demonstrated that variants of ILS, TS, and SA exhibited superior effectiveness and computational efficiency.

Guan et al. (2025) introduce a non-dominated sorting Genetic Algorithm to solve multi-objective FJSPs with batch processing options (FJSBSP), machine speed levels, and transportation that arise in manufacturing settings. The objective is energy efficient production and minimum makespan.

## 3. Methods

In this section, we outline the implementation of the ISA and QL methods, with particular emphasis on methodological developments and other key enhancements. The architecture of our approach is summarised in Figure 1 with details of class attributes and functionality listed in Appendix A.

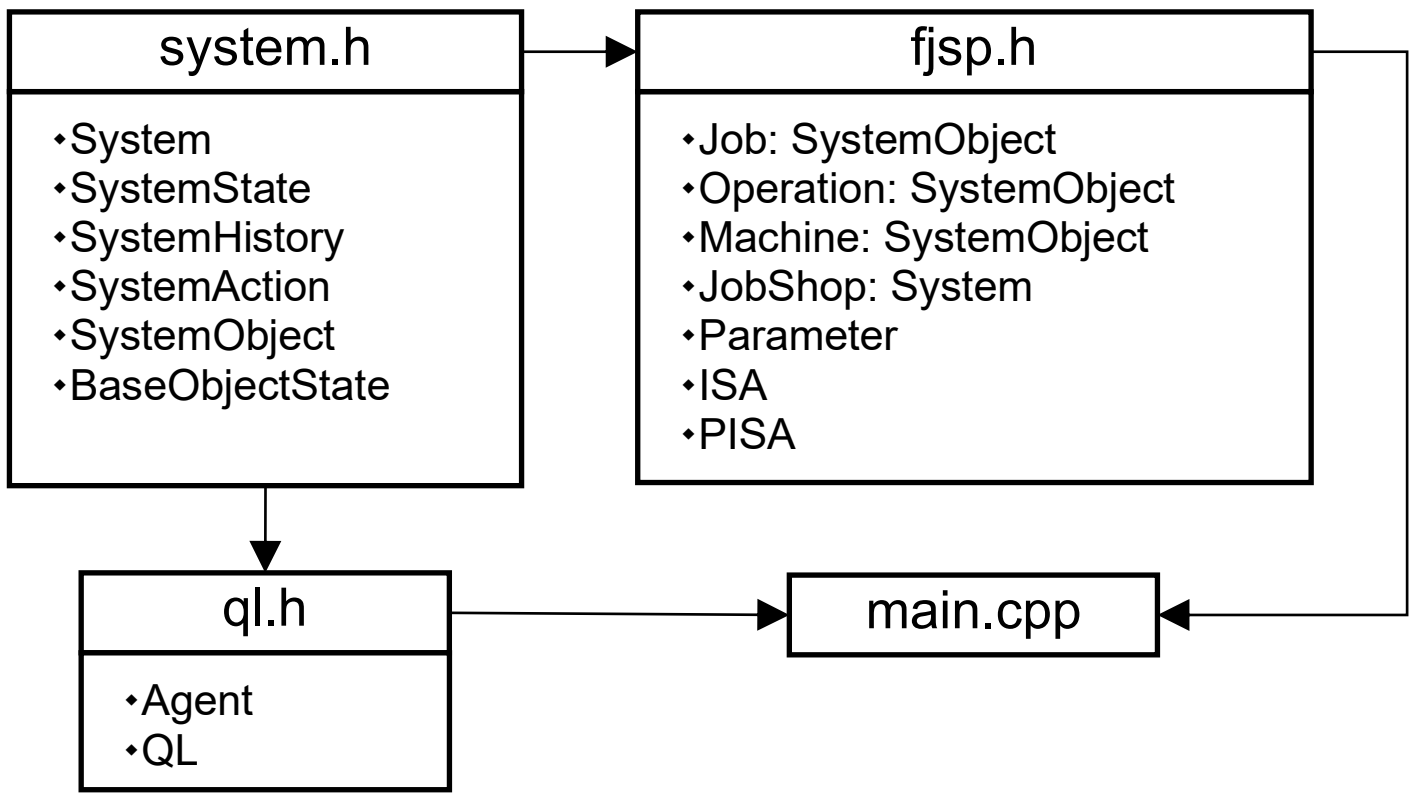

**Figure 1.** Approach architecture

Figure 1 illustrates the organisation of the implementation into three primary modules, each serving a distinct role: system modelling, job shop scheduling, and reinforcement learning through Q-learning. These components are integrated through a central execution point in *main.cpp*.

The first module establishes a foundational framework for tracking and managing system dynamics, historical information, and object states:

- System: Core class that manages the runtime environment and orchestrates the interaction between system components.
- SystemState: Represents the current state of the system at any given point.
- SystemHistory: Tracks historical data of system states or actions, enabling analysis or rollback.
- SystemAction: Encapsulates an action that can be applied to the system.
- SystemObject: Abstract representation of an object within the system.

The second module translates system-level state and action modelling into a concrete application domain: scheduling within a flexible job shop. It also provides the core logic required for scheduling optimization:

- Job: Represents a unit of work consisting of multiple operations.
- Operation: Individual task within a job, typically requiring a specific machine.

- Machine: Represents a resource capable of processing operations.
- JobShop: Manages the overall scheduling environment, combining jobs and machines into a solvable scheduling problem.
- ISA: Implements the Simulated Annealing algorithm for optimizing job scheduling.

The third module applies QL to support decision-making within the job shop environment and comprises the following classes:

- Agent: Learns and makes scheduling decisions using Q-learning.
- QL: Core QL logic, including state/action value updates and policy selection.

It interacts with System and JobShop modules to train agents for better scheduling.

### 3.1. Q-Learning

This subsection outlines the core concepts of Q-learning and their implementation in our framework, with particular emphasis on the careful selection of features identified as most effective in the literature. We adopt a continuous-time modelling approach in which events may occur at any point along the timeline, rather than at fixed discrete intervals. Consequently, no explicit time index is included in the recursive update equation in (1).

Q-learning is fundamentally a feedback-driven reinforcement learning technique (Powell 2019), in which an agent seeks to maximise cumulative reward through interaction with an environment. Each action generates feedback in the form of a reward signal. Positive rewards reinforce desirable behaviour, while negative rewards discourage suboptimal decisions.

The key parameters in QL include the learning rate ($\alpha$), exploration rate ($\varepsilon$), discount factor ($\gamma$), and the number of episodes ($eps$), where $\alpha, \gamma, \varepsilon \in [0,1]$. These parameters govern the learning dynamics and are central to the following recursive update equation:

$$Q(s,a) \;=\; (1-\alpha)Q(s,a) \;+\; \alpha\left(R_{s,a} \;+\; \gamma . \max_{a'}\big(Q(s',a')\big)\right) \qquad (1)$$

In this expression, $Q(s,a)$ denotes the estimated value of action $a$ from state $s$. The term $R_{s,a}$ represents the immediate reward obtained after executing action $a$ from state $s$, while $s'$ denotes the resulting successor state. The expression $\max_{a'}\big(Q(s',a')\big)$ captures the maximum estimated future reward achievable from state $s'$.

Both $s$ and $s'$are associated with a time stamp, and in a continuous-time setting the time corresponding to $s'$ may be equal to that of $s$, depending on the event structure. The discount factor $\gamma$ governs the relative weighting of future rewards. A value of $\gamma = 1$ places full emphasis on long-term returns, whereas $\gamma = 0$ reduces the algorithm to a purely myopic decision rule based only on immediate rewards.

The learning rate $\alpha$ controls the extent to which new information overrides existing estimates. Setting $\alpha = 1$ causes the algorithm to replace the current Q-value with the newly computed estimate, effectively ignoring past experience, whereas $\alpha = 0$ results in no update at all.

The exploration rate $\varepsilon$ determines the balance between exploration and exploitation. When $\varepsilon = 1$, actions are selected uniformly at random, while decreasing values progressively bias the selection toward actions with higher estimated value based on prior experience. Although $\varepsilon$ does not appear explicitly in equation (1), it is incorporated through the ε-greedy policy used to select an action $a_s \in A(s)$, where $A(s)$ denotes the set of feasible actions in state $s$.

$$a_s = \begin{cases} \text{Random choice from } A(s) & \text{with probability } \varepsilon \\ \arg\max_{a \in A(s)} Q(s,a) & \text{with probability } 1-\varepsilon \end{cases} \qquad (2)$$

Conceptually, the ε-greedy policy has limitations in how it structures exploration. In particular, it does not explicitly distinguish between untried and previously sampled actions. A more informative exploration strategy would naturally account for this distinction by partitioning the action space into unseen and seen actions, and by using visitation frequency to guide selection. This would allow the agent to balance exploration of new actions with refinement of value estimates for actions that have already been evaluated.

In this sense, improving or replacing the ε-greedy policy could help mitigate repeated selection of the same actions and promote more structured exploration of the action space. Alternative strategies proposed in the literature include upper confidence bound (UCB) methods and count-based exploration (CBE), both of which explicitly incorporate uncertainty or visitation information into action selection. These are as follows:

$$\text{UCB: } \rho_{a,s} = Q(s,a) - c.\sqrt{\ln(\overline{N}_s)\,/N_{s,a}} \text{ if } N_{n,a} \geq 1, \text{ and } \infty \text{ if } N_{s,a} = 0 \tag{3}$$

$$\text{CBE: } \rho_{a,s} = Q(s,a) - c/\sqrt{N_{s,a} + 1} \tag{4}$$

In (3), $\overline{N}_s$ is the number of visits to state $s$ , i.e., the total number of times any action has been taken in state $s$. In (3) and (4), $\rho_{a,s}$ is the priority of taking a specific action $a$ from state $s$. The action taken at state $s$ is hence computed as follows: $a_s = \arg\max_{a \in A_s} \rho_{a,s}$.

In both UCB and CBE, the priority level for selecting actions is governed jointly by the exploration parameter $c$ and the action-state count $N_{s,a}$ which is continuously updated and records the number of times action $a$ is taken from state $s$. Typically, as $N_{s,a}$ increases the priority decreases. The parameter $c$ directly scales the contribution of the exploration term, functioning similarly to $\varepsilon$. However, unlike $\varepsilon$, it is not constrained to the range [0,1]. A higher value of $c$ encourages greater exploration, but if set too high, the agent may explore excessively, regardless of how often an action has already been evaluated. Conversely, a very low $c$ causes the agent to rely almost entirely on the Q-value, neglecting exploration. In the literature, the range $c \in [1,5]$ is reported as a reasonable balance.

Despite their utility, both UCB and CBE rely on a fixed exploration parameter that can be difficult to tune effectively in practice. Their performance is sensitive to the choice of the exploration constant $c$, and selecting an appropriate value is often problem dependent. In addition, their adaptability to changing problem characteristics or the evolving learning state of the agent is limited, as exploration and exploitation are typically combined through a fixed functional form, imposing a rigid trade-off that may not be optimal across all stages of learning (Zandvakili and Mansouri, 2025)

To address these limitations, a more flexible, instance-aware mechanism is desirable, capable of dynamically balancing exploration and exploitation based on observed data. In this context, both the action-value estimate $Q(s,a)$and the visitation count $N_{s,a}$ can be jointly utilised within a Pareto-based selection (PBS) framework (Van Moffaert and Nowe (2014), Drugan and Nowe (2013)). This approach promotes actions with high estimated value while simultaneously encouraging exploration of less-visited state–action pairs, thereby supporting more balanced learning.

Within this framework, candidate actions are identified as those that are non-dominated in a two-dimensional objective space defined by $(Q(s,a), F(N_{s,a}))$, where $F(\cdot)$ is a monotonically decreasing function of visitation frequency (e.g., $F(N_{s,a}) = -N_{s,a}$or $F(N_{s,a}) = 1/(N_{s,a} + 1)$). The ideal (utopia) point corresponds to high $Q(s,a)$and low $N_{s,a}$. An action is considered non-dominated if no other action simultaneously achieves a higher value estimate and a lower visitation count (see Figure 2). Once the non-dominated set is constructed, an action is selected uniformly at random from this set.

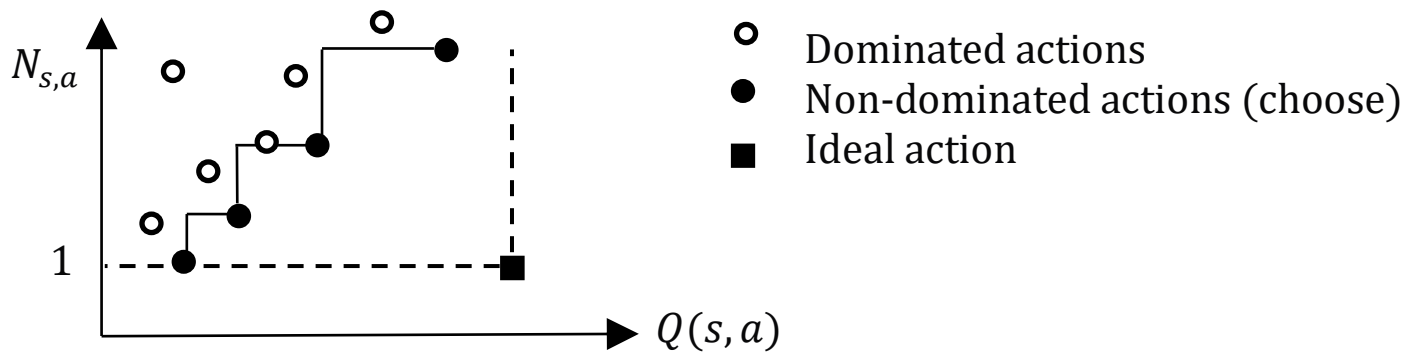

**Figure 2.** Demonstration of PBS

**Stochastic Rewards**. In stochastic FJSP settings, rewards are non-deterministic, and the same action, when applied from a given state, can yield different outcomes across episodes. Consequently, estimating both the expected value $\mathbb{E}[Q_n(s,a)]$ and the variance $\mathbb{V}[Q_n(s,a)]$ provides a more robust basis for decision-making under uncertainty. Let us denote $Q_n^1(s,a)$ as the expected value $\mathbb{E}[Q_n(s,a)]$, and $Q_n^2(s,a)$ as its variance $\mathbb{V}[Q_n(s,a)]$ after $n$ samples. They may be updated incrementally as follows, for $n \geq 2$:

$$Q_n^1(s,a) = Q_{n-1}^1(s,a) + (1/n)(Q_n^{\text{new}} - Q_{n-1}^1(s,a)) \tag{5}$$

$$Q_n^2(s,a) = (1/n)\left[Q_{n-1}^2(s,a) + \left(Q_n^{\text{new}} - Q_{n-1}^1(s,a)\right)\left(Q_n^{\text{new}} - Q_n^1(s,a)\right)\right] \tag{6}$$

Where, $n = N(s,a)$ and $Q_n^{\text{new}} = R_{s,a,n} + \gamma . \max_{a'}\left(Q_{n-1}^1(s',a')\right)$, $Q_1^1(s,a) = Q_1^{\text{new}}$ and $Q_1^2(s,a) = 0$. Notably, equation (6) is a variation of Welfords' approach (Welford 1962).

As more data is collected during learning, the estimate $\mathbb{E}[Q(s,a)]$ should stabilize, while the variance $\mathbb{V}[Q(s,a)]$ should typically decrease. In QL, the guiding principle is to select actions that are either clearly good or have the potential to be good. A low expected value coupled with low variance suggests high confidence in a poor decision. Generally, such actions should be avoided. Conversely, a high expected value with low variance indicates a reliable, favourable choice. High variance, regardless of the expected value, signals uncertainty and may warrant continued exploration of that action.

Given the previously discussed limitations of the e-greedy approach, it is natural to adopt the following priority metric for comparing and selecting actions:

$$\rho_{a,s} = \mathbb{E}[Q(s,a)] + z.\sqrt{\mathbb{V}[Q_n(s,a)]} \tag{7}$$

This approach balances exploitation and uncertainty-guided exploration, where the parameter $z$ controls the degree of exploration. This idea naturally extends to a Pareto-based selection approach where actions are evaluated in a two-dimensional objective space defined by $(\mathbb{E}[Q(s,a)], -\mathbb{V}[Q(s,a)])$. This Pareto front consists of actions that are non-dominated, those for which no other action has both higher expected value and lower variance. The utopia point in this space corresponds to actions with high expected returns and low uncertainty. One action is then selected at random from this set of Pareto-optimal candidates. Notably, this approach implicitly incorporates the visitation count $N_{s,a}$ through the variance, as uncertainty naturally decreases with increased sampling.

**Recording Q-Values**. In Q-learning, each state–action pair is associated with a corresponding Q-value. However, in the FJSP, as in many combinatorial decision problems, the state space is extremely large, making it impractical to predefine and store all values in a static Q-table. In our C++ implementation, we therefore use an std::unordered_map<string, vector<double>> to store Q-values dynamically, denoted as $Q_m(s,a)$ or $Q_{m,n}^1(s,a)$ and $Q_{m,n}^2(s,a)$, depending on the formulation. In C++, unordered_map is a hash table–based associative container that stores key–value pairs and provides average-case $O(1)$ complexity for lookup, insertion, and deletion operations. This structure effectively acts as a dictionary that is updated as new states are encountered during learning. Each state entry is associated with a fixed-length vector, where the dimension corresponds to the total number of available actions.

To uniquely identify each system state, we define a custom, application-specific hash function that generates a string-based key encoding the complete status of the system at a given time. This string must compactly and unambiguously represent the current time and the status of all jobs in the system. The resulting string consists of the current time followed by a status substring for each job, with entries separated by the delimiter "|". Each job substring encodes its current state. If a job is completed, it is represented as "c". If a job is in progress and the next operation has not yet started, the substring is "O#,p", where O# denotes the operation index. If the operation is currently being processed, the substring is "O#,M#,rempt", indicating the operation index, the assigned machine, and the remaining processing time (e.g., "1,1,14"). The terminal state, in which all jobs are completed, is represented as: t|c|c|...|c. The following code, (17|2,2,4|1,1,14|3,p|2,3,1), illustrates these conventions. This string encodes the system state at time 17, where:

- J1_O2 is running on machine M2 with 4 units of processing time remaining.
- J2_O1 is running on M1 with 14 units remaining.
- J3_O3 is pending.
- J4_O2 is running on M3 with 1 unit remaining.

From this, we can infer that M1, M2, and M3 are currently busy. M3 will be free at time 18 (17+1), and M1 at time 31 (17+14). M2 will be free at time 21 (17+4).

**Actions.** Because actions are inherently application-specific, we introduce a generic Action class. This class encapsulates the action itself using a std::any attribute, alongside its associated reward and timestamp. By abstracting the action in this way, the reinforcement learning algorithm can operate on actions without requiring knowledge of their internal structure. Any application-specific interpretation or processing is deferred to the appropriate context (e.g., the FJSP class), allowing the core RL logic to remain modular, general, and reusable.

In the context of the FJSP, an action corresponds to assigning a particular operation to a specific machine. The set of feasible assignments is known in advance. The total number of feasible assignments defines the size of the action space, and consequently, the number of entries associated with each state in the Q-table.

**System and States**. We implement generic System and SystemState classes to represent and manage the state of a system. The System class maintains both an initial and a current state and primarily serves as a base abstraction for defining application-specific behaviour. It also provides the core interfaces required for integration with Q-learning algorithms, ensuring that domain logic remains decoupled from the learning framework.

The SystemState class captures the state of all entities within the system. Internally, it stores a dictionary defined as unordered_map<string, StateRecord> records, where each entry maps a unique identifier to a corresponding StateRecord. Each StateRecord contains a pointer to the associated object and a variable representing its current state. To support extensibility, it also includes a more_info string field, which can store additional object-specific metadata or contextual details as needed.

In the context of the FJSP, the system models machines, jobs, and their associated operations, each represented as a derived Object type. Jobs and operations can exist in one of three states: pending, in progress, or complete, while machines are either idle or busy. Supplementary information—such as job or operation completion times, or the remaining processing time of an operation—is stored in more_info. For machines, this field captures the time until the resource becomes available.

**Transitioning Between States.** The System class defines a set of pure virtual functions that establish the interface for core operations, including state transitions, action execution and evaluation, decision point identification, state hashing, action logging, and reward computation. No default implementation is provided; instead, these functions must be implemented by the user to reflect the specific logic and constraints of the target application.
The transition function determines the next system state resulting from a given action. In the context of the FJSP, this involves iterating over all operations and updating their status relative to the current time. Operations that were in progress and have completed by the specified time are marked as complete, and their assigned machines are set to idle. If an operation is still in progress, its state remains unchanged, but its remaining processing time is updated. Completed operations require no further modification.

The take_action function assigns a selected operation to a specific machine and updates the corresponding states of both entities. The resulting system time advances to reflect the completion of the action and the subsequent state change, ensuring consistency between action execution and system evolution.

**Identifying the Next Decision Point.** This requirement is particularly important in continuous-time models. After an action is taken, the next decision point does not necessarily advance in time—it may coincide with the current time. Conversely, simply jumping to the next event can be inappropriate, as that event may not correspond to a point where a new decision is feasible. It is therefore necessary to explicitly identify the next time at which a valid action can be taken.

In the context of the FJSP, this is achieved by examining the state records of all machines within the current SystemState. For each busy machine, the remaining processing time—stored in the more_info field—is used to compute its next availability time. For idle machines, the availability time is simply the current time. The next decision point is then defined as the minimum of these availability times across all machines.

**Rewards.** After an action is executed, an associated reward must be assigned. However, the true value of an action, such as assigning an operation to a specific machine at a given time, may not be immediately meaningful or informative, and its impact may only become apparent later in the schedule.

Several basic reward formulations have been proposed in the literature. For example, one approach defines the reward as $R(s,a) = -C_{\max}(s')$, where the reward corresponds to the makespan of the resulting (partial) schedule after the action is applied. An alternative, $R(s,a) = C_{\max}(s) - C_{\max}(s')$, captures the incremental change in makespan induced by the action. In both cases, rewards are negative to reflect the objective of minimising makespan. These options are inherently myopic, as they focus only on immediate outcomes and do not capture how a current decision influences future scheduling flexibility or downstream performance. A contrasting approach is to use a sparse reward structure, where $R(s,a) = 0$ for all non-terminal transitions, and $R(s,a) = -C_{\max}(s')$ only when a terminal state is reached. Under this formulation, the reward signal is deferred until the completion of the schedule. While this avoids short-sighted feedback, it introduces a challenging credit assignment problem, as it becomes difficult to attribute the outcome to individual actions. It can also hinder efficient exploration, since the learning signal is infrequent and delayed.

Bengio et al., (2021) suggests that it may be beneficial to design a surrogate reward signal that encourages intermediate accomplishments. This approach, known as *reward shaping*, is an active area of research in reinforcement learning. A key requirement of reward shaping is that it must preserve the optimal policy and avoid introducing bias that could lead the agent toward suboptimal behaviour. For the FJSP, we identify several forms from the literature worth further investigation:

$$R(s,a) = -Cmax(s') + \lambda.\mathrm{f}(s,a) \tag{8}$$
$$R(s,a,s') = R(s,a) + \phi(s') - \phi(s) \tag{9}$$

In Equation (8), the reward is augmented with an auxiliary shaping term $\mathrm{f}(s,a)$, weighted by a factor $\lambda$, which incorporates domain-specific signals such as partial progress, resource utilisation, or idle time reduction. This provides additional intermediate feedback while still anchoring learning to the primary objective of makespan minimisation. In Equation (9), $\phi(s)$ provides a score for the given state and the difference $\phi(s') - \phi(s)$ is added to the reward to provide additional feedback based on the transition between states. To guide learning in the FJSP environment, $\mathrm{f}(s,a)$ can be intuitively defined in several meaningful ways:

$$\mathrm{f}(s,a) = O_s^{\mathrm{done}}/O^{\mathrm{total}} \text{ where } O^{\mathrm{total}} = \textstyle\sum_{j\in J}|\mathcal{O}_j| \tag{10}$$

$$\mathrm{f}(s,a) = CPT/TPT \text{ or } \mathrm{f}(s,a) = -RPT/TPT \text{ where } TPT = \textstyle\sum_{j\in J}\sum_{k\in\mathcal{O}_j} p_{j,k} \tag{11}$$

$$\mathrm{f}(s,a) = -I(s') \tag{12}$$

The first option defines the progress ratio as the number of completed operations ($O^{\mathrm{done}}$) divided by the total number of operations $O^{\mathrm{total}}$. However, partially complete operations do not contribute to that measure. To address this limitation, we define $O_s^{\mathrm{done}}$ as $\sum_{j\in J}\sum_{k\in\mathcal{O}_j} prop_{j,k}^s$ where $prop_{j,k}^s \in [0,1]$ denotes the degree of completion (a.k.a. the work-in-progress proportion)of $o_{j,k}$ in state $s$.

The second option defines progress in terms of time, for example as the ratio of completed processing time (CPT) to total processing time (TPT), or equivalently through the remaining processing time (RPT) relative to TPT. However, for the FJSP, these measures are problematic because processing times are machine dependent. As a result, both CPT and RPT cannot be computed reliably in advance or in a machine-independent manner. The third shaping signal considers the negative of idle time within the partial schedule, thereby encouraging decisions that maximise machine utilisation. Idle time typically arises because of suboptimal sequencing and assignment decisions and reducing it can lead to more efficient schedules.

**Extracting the Optimal Schedule.** According to the existing literature, an optimal policy can be obtained by starting from the initial state and iteratively selecting the action with the highest Q-value until a terminal state is reached. However, our numerical experiments indicate that this procedure does not necessarily recover the best-performing solution (i.e., schedule) encountered during learning. This observation is based on an analysis of the Q-table, which revealed terminal states with superior performance that were not obtained through standard greedy policy extraction. This limitation is particularly evident in deterministic settings, where the learned value function may still fail to consistently encode the globally best trajectory.

In retrospect, several strategies can be used to address this issue and identify the true best solution discovered during training. The most effective approach adopted in our implementation is to maintain a record of the agent's actions throughout each episode. As actions are executed, they are appended to a sequence representing the current schedule. Once a terminal state is reached, this sequence defines a complete solution, which is then evaluated and compared against the best solution found so far. If it yields a higher performance score, it replaces the current best solution. At the end of each episode, the action list is reset. This approach is used in the numerical experiments reported later. Alternatively, one could construct a state transition tree (Q-tree), as illustrated in Figure 3. This structure is built from transition data collected during learning, starting from a root node corresponding to the initial state and progressively expanding as new state–action transitions are observed. The tree ultimately branches into terminal states, enabling systematic reconstruction and evaluation of complete solution paths.

While the Q-tree approach provides useful insight into the structure of the learned policy and the exploration dynamics, it is primarily of academic interest. It was not adopted in our implementation due to scalability concerns, as our focus is on computationally efficient methods suitable for larger problem instances.

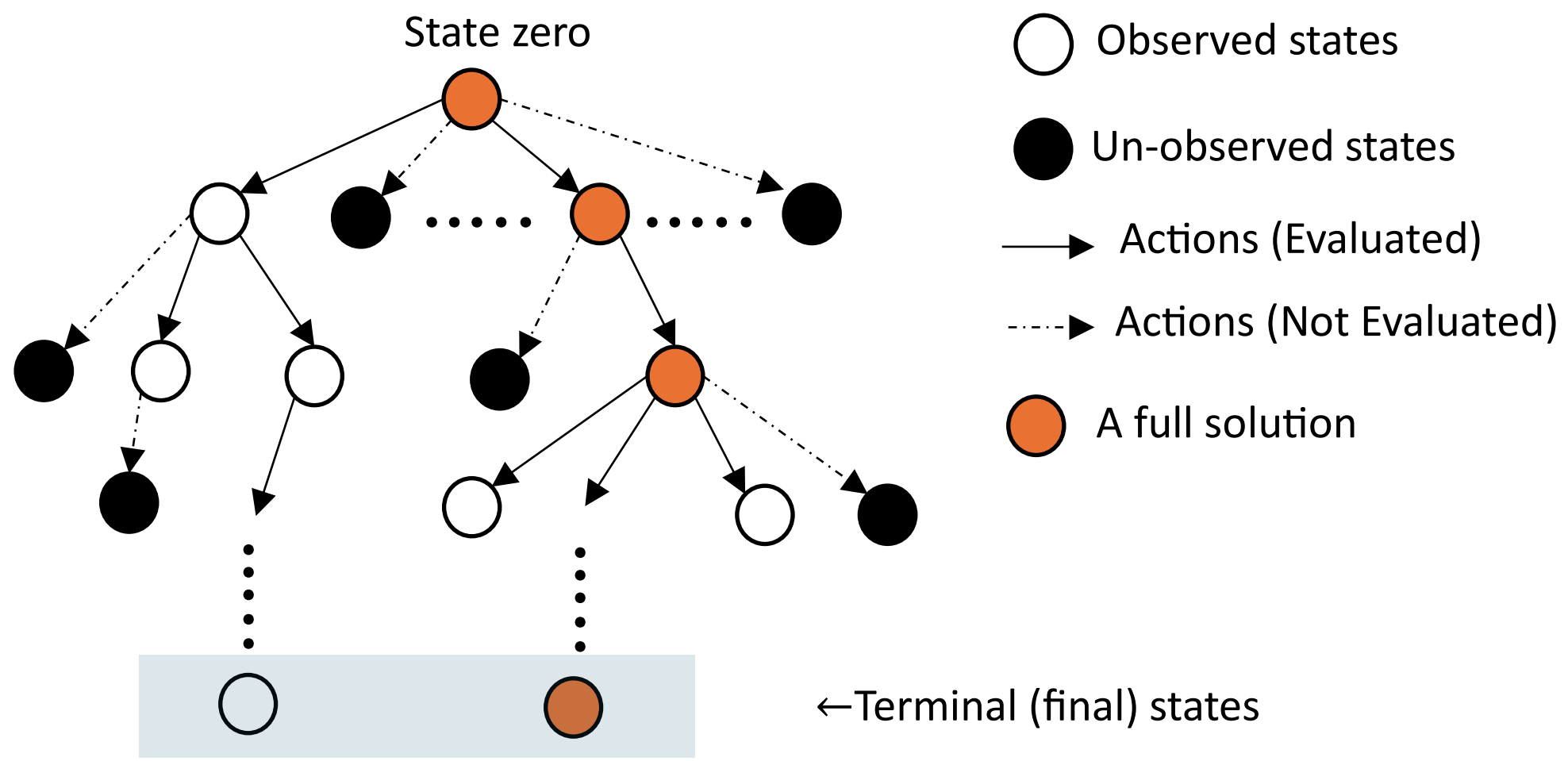

**Figure 3.** Demonstrative state transition tree

**3.2. Simulated Annealing with Intensification**

Our C++ implementation of ISA, referred to as PISA, is a parallel (i.e., distributed) variant designed to exploit the computational capabilities of modern multi-processor systems. The central idea behind PISA is to extend the search capacity of a traditional serial ISA by distributing the computational workload across multiple solvers. In each iteration, multiple ISA solvers independently explore different regions of the solution space. At the end of each temperature step, all solvers report the best solution they have identified, allowing this information to be shared and used to guide subsequent search progress. This cooperative mechanism enables more effective exploration while retaining the core structure of the underlying ISA framework.

A central feature of our approach is the use of a priority ordering (PO) solution representation (Burdett et al., 2019). A priority order is defined as a finite sequence of pairs $(o, m)$, where $o \in \mathcal{O}$ denotes an operation and $m \in \mathcal{M}$ denotes a machine. The sequence length equals the total number of operations to be scheduled, i.e., $|\mathcal{O}|$. Each pair must satisfy machine eligibility constraints, meaning that the selected machine must be capable of processing the corresponding operation. In addition, the sequence must preserve the precedence relationships associated with each job. Specifically, for any $o, o' \in \mathcal{J}_i$, if $o \preccurlyeq o'$ then $(o, m_o) \preccurlyeq (o', m_{o'})$ within the priority sequence. Let $\mathcal{S} \subseteq \mathcal{O} \times \mathcal{M}$ denote the set of all eligible operation-machine pairs, defined as $\mathcal{S} = \{(o, m) \in \mathcal{O} \times \mathcal{M} | o \text{ can be processed on } m\}$. Then the priority order is an element of the set $\mathcal{S}^{|\mathcal{O}|}$. In C++, this structure is naturally represented as std:: vector<pair<Operation*, Machine*>>. This enables us to address the scheduling problem without relying on a disjunctive graph (DG) and its associated complexity, providing a more lightweight and computationally efficient alternative.

Traditionally, a disjunctive graph is used to manage operational scheduling and detect precedence violations (e.g., cycles), typically by transforming machine sequences into disjunctive arcs and applying graph-based algorithms such as longest-path or Bellman–Ford procedures. In contrast, our priority ordering directly induces the required topological structure. Consequently, scheduling can be implemented through a simple dispatching policy that incrementally updates machine availability as operations are processed.

Schedule evaluation in the stochastic setting differs fundamentally from the deterministic case. The most robust approach involves simulating scenarios until both the mean and variance of the makespan converge, ensuring that sufficient stochastic variability is captured. In each scenario, operation processing times are sampled from their respective probability distributions.

However, this convergence-based method is inherently sequential and computationally expensive. To improve scalability, a parallel alternative can be employed in which a fixed number of scenarios is used to balance accuracy with computational efficiency. In practice, a preliminary convergence study can be performed to estimate an appropriate number of scenarios, which is then used to parameterise the parallel evaluation scheme.

Perturbation of the current solution is a fundamental component of SA, and the design of the perturbation mechanism plays a critical role in overall algorithmic performance. After careful evaluation, two strategies were identified as particularly effective: shift operations and machine reassignments, both of which are complemented by local search to enhance solution quality and convergence behaviour. This coupling is a hallmark of SA with intensification.

The shift operator relocates a selected operation either to the left or right within the priority sequence. The magnitude of the shift is chosen randomly but is constrained by the positions of the operation's immediate predecessor and successor within the same job, ensuring feasibility with respect to precedence relations. The machine reassignment operator replaces the currently assigned machine with an alternative feasible machine, provided that such an assignment exists.

Our perturbation operators are compound moves in the sense that they span multiple neighbourhood structures. For instance, after an operation is shifted within the priority sequence, its existing machine assignment may no longer be appropriate or may become suboptimal in the new context. Applying a local search to revise the machine assignment is therefore straightforward and helps eliminate clearly inferior local choices. Similarly, after a machine reassignment, the current position of the operation–machine pairing may be suboptimal. In such cases, a local search over feasible positions is beneficial.

To further improve convergence, local search is applied once per temperature step, focusing on the current best solution. Two local search strategies are used:

1. For a randomly selected operation, all feasible shift positions within the sequence are evaluated.
2. For a given operation, all valid machine assignments are explored.

In both cases, the best-improving move is identified and applied if it yields an improvement. To manage computational cost, only one of the two local search strategies is selected at random and applied at each temperature step.

Our analysis supports the use of the traditional SA acceptance criterion, in which newly generated solutions are evaluated and accepted or rejected according to their relative quality. An alternative strategy, threshold acceptance, has also proven effective in certain settings. Furthermore, advanced perturbation mechanisms that focus exclusively on critical operations have been investigated. While such targeted strategies can improve performance by concentrating search effort on influential regions of the solution space, they may also degrade overall solution quality by reducing diversity and constraining broader exploration.

## 4. Results

This section presents a series of numerical experiments designed to demonstrate the effectiveness of the proposed methods. The evaluation has two primary objectives: (i) to assess whether the streamlined parallel, collaborative simulated annealing (SA) implementation performs competitively against established benchmark approaches, and (ii) clarify how the Q-learning approach compares with both PISA and the benchmark methods.

As benchmark instances, we employ the Dauzère-Pérès and Paulli FJSP dataset (Dauzere-Peres and Paulli, 1997), which is widely used in the literature and provides a robust basis for comparative analysis. These instances are particularly challenging, as each job comprises more operations than there are machines, and each operation can be processed on multiple alternative machines. This combination of high routing flexibility and complex precedence structure makes them well suited for evaluating advanced scheduling algorithms. Given their difficulty and broad representation of key problem characteristics, we consider performance on these instances to be indicative of overall algorithmic effectiveness.

Preliminary experiments were conducted on a local workstation equipped with an 11th Generation Intel Core i7 processor (base frequency 3.3 GHz, turbo up to 4.8 GHz), featuring four physical cores (eight logical threads) and 16 GB of RAM, running Windows 11. Subsequent large-scale experiments were performed in a high-performance computing (HPC) environment using AMD EPYC 9684X processors, each providing 96 physical cores with a base clock speed of 2.55 GHz and a maximum boost frequency of 3.7 GHz.

### 4.1. FJSP with Deterministic Processing Times

Details of the FJSP instances are displayed in Table 3. Key summary statistics of the 18 instances are shown, accompanied by best solutions reported by the two bespoke algorithms of theirs, denoted CMJS2 and CBest. The later (superior) one is a Tabu Search algorithm.

**Table 3.** Problem instances

| Inst. | #Jobs | #Mach | #Opns | #Options | Ratio: #Options To #Opns | LB | CMJS2 | Cbest | cpu (sec) |
|---|---|---|---|---|---|---|---|---|---|
| **1** | 10 | 5 | 196 | 221 | 1.13 | 2505 | 2642 | 2530 | 260.1 |
| **2** | 10 | 5 | 196 | 332 | 1.69 | 2228 | 2317 | 2244 | 251.2 |
| **3** | 10 | 5 | 196 | 501 | 2.56 | 2228 | 2263 | 2235 | 256.2 |
| **4** | 10 | 5 | 196 | 221 | 1.13 | 2503 | 2648 | 2565 | 263.9 |
| **5** | 10 | 5 | 196 | 332 | 1.69 | 2189 | 2319 | 2229 | 317.9 |
| **6** | 10 | 5 | 196 | 501 | 2.56 | 2162 | 2263 | 2216 | 3778.2 |
| **7** | 15 | 8 | 293 | 364 | 1.24 | 2187 | 2473 | 2408 | 305.7 |
| **8** | 15 | 8 | 293 | 708 | 2.42 | 2061 | 2150 | 2093 | 260.5 |
| **9** | 15 | 8 | 293 | 1182 | 4.03 | 2061 | 2127 | 2074 | 685.2 |
| **10** | 15 | 8 | 293 | 364 | 1.24 | 2178 | 2460 | 2362 | 251.6 |
| **11** | 15 | 8 | 293 | 708 | 2.42 | 2017 | 2124 | 2078 | 894.7 |
| **12** | 15 | 8 | 293 | 1182 | 4.03 | 1969 | 2073 | 2047 | 839.6 |
| **13** | 20 | 10 | 387 | 518 | 1.34 | 2161 | 2360 | 2302 | 328.9 |
| **14** | 20 | 10 | 387 | 1156 | 2.99 | 2161 | 2224 | 2183 | 681.6 |
| **15** | 20 | 10 | 387 | 1941 | 5.02 | 2161 | 2295 | 2171 | 1732.4 |
| **16** | 20 | 10 | 387 | 518 | 1.34 | 2148 | 2391 | 2301 | 341.7 |
| **17** | 20 | 10 | 387 | 1156 | 2.99 | 2088 | 2211 | 2169 | 491.8 |
| **18** | 20 | 10 | 387 | 1941 | 5.02 | 2057 | 2160 | 2139 | 1944.7 |

In the preliminary exploratory phase, PISA was systematically evaluated under a range of parameter settings to assess its behaviour and performance. Two configurations are particularly noteworthy. The first corresponds to shorter runs with a stronger greedy bias ($t_{\mathrm{I}} = 1$, $t_{\mathrm{F}} = 0.001$, $maxSteps = 100$, $maxIter = 5000$). The second involves longer runs with higher initial temperature and increased perturbation intensity ($t_{\mathrm{I}} = 10$, $t_{\mathrm{F}} = 0.001$, $maxSteps = 100$, $maxIter = 10000$). The parameter maxIter, which controls the number of perturbations performed at each temperature level, is set relatively high in both cases. This reflects the difficulty of the Dauzère-Pérès and Paulli benchmark instances, which require extensive exploration to achieve high-quality solutions

Within PISA, four parallel ISA searches were executed to achieve effective utilisation of the available computational resources while limiting communication and synchronisation overhead. Although a larger number of parallel instances could potentially yield further performance improvements, this was intentionally avoided to ensure a fair comparison with the Q-learning approach, which does not benefit from a comparable degree of parallelism.

The aggregated results are reported in Tables B1 and B2 in Appendix B. The percentage difference (%DIFF) is computed as $100 \times (\mathrm{PISA} - \mathrm{CBest})/\mathrm{CBest}$and $100 \times (\mathrm{PISA} - \mathrm{CMSJ2})/\mathrm{CMSJ2}$, respectively, with negative values indicating superior performance relative to the benchmark.

Figures 4 and 5 provide a visual comparison of performance relative to the benchmark solutions, while Figure 6 illustrates the consistency of the solutions produced by PISA. Differences in runtime between the "short" and "long" configurations are presented in Figure 7. We next evaluated the Q-learning approach under a range of hyperparameter settings informed by the literature. The learning rate and exploration parameter ($\epsilon$) were varied within the interval [1, 0.1] using an exponential decay schedule, while the discount factor was explored over the range [0.1, 1].

For action selection, both $\epsilon$-greedy and Pareto-based strategies were employed. Upper Confidence Bound and CBE-based approaches were excluded due to their known limitations in this context, as well as to maintain focus and clarity in the experimental design.

With respect to the reward signal, multiple formulations were considered, including incremental changes in makespan, incumbent makespan values, and the ratio of completed operations. Aggregated results, based on ten independent trials and 20,000 learning episodes, are reported in Appendix B (Tables B3–B6). Performance is evaluated using the percentage difference metric, computed as 100×(QL – PISA)/PISA).

Early-stage experiments were conducted using 1,000, 5,000, and 10,000 training episodes to assess scalability. As Q-learning performance improved consistently with increasing numbers of episodes, we selected 20,000 episodes—the maximum supported by the memory constraints of the local workstation—for subsequent experiments.

The corresponding results are presented in Tables B3 and B4. The "%PISA" column reports the relative performance difference with respect to PISA (best), with negative values indicating superior performance. Notably, the Pareto-based selection strategy consistently outperformed the $\epsilon$-greedy approach across all instances, as reflected by the negative values in the "%E-GRDY" column. We also observed that variations in the reward function had minimal impact on overall performance. Consequently, results are reported using a single representative reward function.

***Analysis***. To evaluate the effectiveness of the proposed methods, formal comparisons were conducted both against the benchmark algorithms and between the methods themselves. Tables 4 and 5 presents pairwise comparisons for both best-case and average-case performance.

**Table 4.** Pairwise comparison and statistical tests for best case

| A | B | #better | | mean diff | stdev | two-sided | | | one-sided | |
|---|---|---|---|---|---|---|---|---|---|---|
| | | A | B | B-A | | t | p | sig | p | sig |
| PISA (short) | PISA (long) | 5 | 13 | -6.056 | 50.820 | 0.506 | 0.620 | × | 0.310 | × |
| QL (egrdy) | QL (pbs) | 1 | 17 | -31.722 | 34.958 | 3.850 | 0.001 | ✓ | 0.001 | ✓ |
| QL (egrdy) | PISA (short) | 9 | 9 | -55.556 | 117.849 | 2.000 | 0.062 | × | 0.031 | ✓ |
| QL (egrdy) | PISA (long) | 8 | 10 | -61.611 | 111.095 | 2.353 | 0.031 | ✓ | 0.015 | ✓ |
| QL (pbs) | PISA (short) | 10 | 8 | -23.833 | 106.755 | 0.947 | 0.357 | × | 0.178 | × |
| QL (pbs) | PISA (long) | 9 | 9 | -29.889 | 92.139 | 1.376 | 0.187 | × | 0.093 | × |
| Cbest | QL (egrdy) | 18 | 0 | 154.722 | 90.252 | 7.273 | 0.000 | ✓ | 0.000 | ✓ |
| Cbest | QL (pbs) | 18 | 0 | 123.000 | 89.172 | 5.852 | 0.000 | ✓ | 0.000 | ✓ |
| Cbest | PISA (short) | 17 | 1 | 99.167 | 65.679 | 6.406 | 0.000 | ✓ | 0.000 | ✓ |
| Cbest | PISA (long) | 17 | 1 | 93.111 | 79.452 | 4.972 | 0.000 | ✓ | 0.000 | ✓ |
| CMJS2 | QL (egrdy) | 17 | 1 | 90.611 | 81.339 | 4.726 | 0.000 | ✓ | 0.000 | ✓ |
| CMJS2 | QL (pbs) | 15 | 3 | 58.889 | 84.304 | 2.964 | 0.009 | ✓ | 0.004 | ✓ |
| CMJS2 | PISA (short) | 11 | 7 | 35.056 | 84.268 | 1.765 | 0.096 | × | 0.048 | ✓ |
| CMJS2 | PISA (long) | 11 | 7 | 29.000 | 90.707 | 1.356 | 0.193 | × | 0.096 | × |

**Table 5.** Pairwise comparison and statistical tests for average case

| A | B | #better | | mean diff | stdev | two-sided | | | one-sided | |
|---|---|---|---|---|---|---|---|---|---|---|
| | | A | B | B-A | | t | p | sig | p | sig |
| PISA (short) | PISA (long) | 9 | 9 | 19.833 | 111.476 | 0.755 | 0.461 | × | 0.230 | × |
| QL (egrdy) | QL (pbs) | 0 | 18 | -36.667 | 24.144 | 6.443 | 0.000 | ✓ | 0.000 | ✓ |
| QL (egrdy) | PISA (short) | 10 | 8 | 40.778 | 152.416 | 1.135 | 0.272 | × | 0.136 | × |
| QL (egrdy) | PISA (long) | 14 | 4 | 60.611 | 132.158 | 1.946 | 0.068 | × | 0.034 | ✓ |
| QL (pbs) | PISA (short) | 11 | 7 | 77.444 | 143.691 | 2.287 | 0.035 | ✓ | 0.018 | ✓ |
| QL (pbs) | PISA (long) | 14 | 4 | 97.278 | 114.606 | 3.601 | 0.002 | ✓ | 0.001 | ✓ |
| Cbest | QL (egrdy) | 18 | 0 | 180.222 | 100.684 | 7.594 | 0.000 | ✓ | 0.000 | ✓ |
| Cbest | QL (pbs) | 18 | 0 | 143.556 | 94.264 | 6.461 | 0.000 | ✓ | 0.000 | ✓ |
| Cbest | PISA (short) | 18 | 0 | 221.000 | 87.691 | 10.692 | 0.000 | ✓ | 0.000 | ✓ |
| Cbest | PISA (long) | 18 | 0 | 240.833 | 126.596 | 8.071 | 0.000 | ✓ | 0.000 | ✓ |
| CMJS2 | QL (egrdy) | 17 | 1 | 116.111 | 91.368 | 5.392 | 0.000 | ✓ | 0.000 | ✓ |
| CMJS2 | QL (pbs) | 16 | 2 | 79.444 | 87.040 | 3.872 | 0.001 | ✓ | 0.001 | ✓ |
| CMJS2 | PISA (short) | 18 | 0 | 156.889 | 106.342 | 6.259 | 0.000 | ✓ | 0.000 | ✓ |
| CMJS2 | PISA (long) | 18 | 0 | 176.722 | 133.280 | 5.626 | 0.000 | ✓ | 0.000 | ✓ |

Paired *t*-tests were used to assess whether the observed differences are statistically significant. The null hypothesis assumes no effect; that is, there is no difference in the mean makespan across replications for each test instance. A *p*-value below 0.05 is taken as evidence to reject the null hypothesis, indicating a statistically significant difference. One-sided tests are employed to evaluate whether Algorithm B yields lower makespan values than Algorithm A.

In summary, both PISA and QL consistently generated competitive solutions, although neither matched the performance of CBest under the current experimental conditions. PISA surpassed CBest on a small number of instances, whereas QL did not. The best solutions obtained by PISA were comparable to those of CMJS2, with no statistically significant difference detected.

PISA (short) and PISA (long) produced broadly similar solution quality, with only minimal statistical distinction between them. Nevertheless, PISA (long) achieved slightly better average performance and delivered superior solutions in thirteen cases. Both PISA and QL (pbs) significantly outperformed QL (ε-greedy), while QL (pbs) achieved performance comparable to both PISA variants. Although PISA retained a slight overall numerical advantage—reflected by average differences of –23.833 and –29.889—the aggregate analysis indicates that these differences were relatively modest. Indeed, QL (ε-greedy) was not significantly worse on average, and QL (pbs) was statistically superior to PISA, with mean improvements of 77.44 and 97.278, respectively.

Q-learning demonstrated unexpectedly strong performance on several instances. Depending on whether ε-greedy exploration or Pareto-based selection was applied, it outperformed PISA on average in eight of the eighteen instances. The number of trials in which Q-learning achieved superior performance is presented in Figure 8, while Figure 9 provides a direct comparison of best makespan values. Figure 10 reports the average number of states stored in the Q-table after 20,000 learning episodes, together with the number of complete solutions that could be extracted, illustrating the substantial size of the underlying state space.

PISA's principal advantage lies in computational speed, as Q-learning was notably slower on some instances. Even so, Q-learning's overall performance remains noteworthy. Its substantially lower standard deviations indicate greater consistency across runs compared with PISA. Although Q-learning incurred higher CPU times, this was expected given the absence of parallelization. The disparity became more pronounced on larger instances. By contrast, PISA benefited from four ISA solvers operating in parallel, which contributed to its runtime advantage.

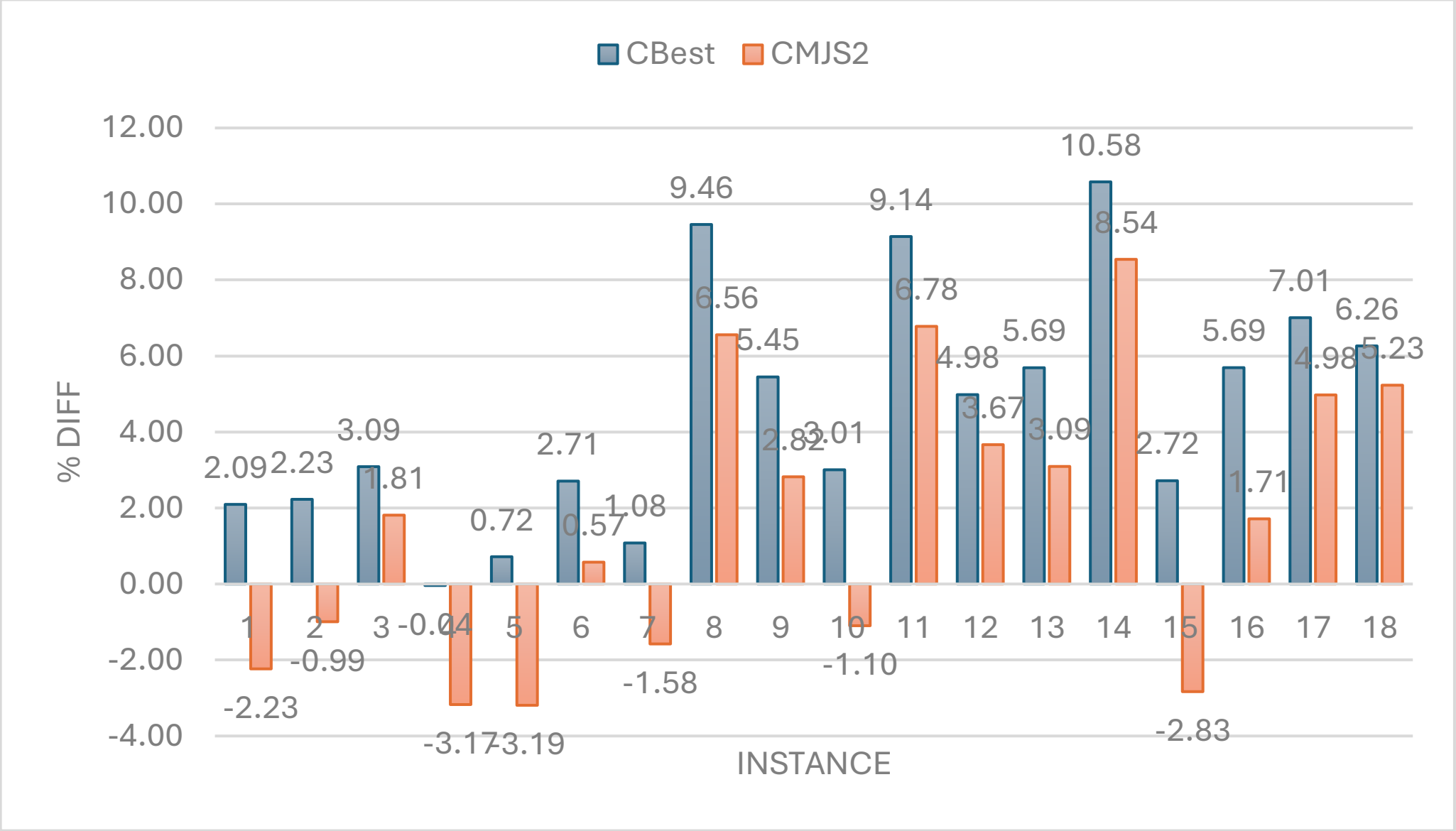

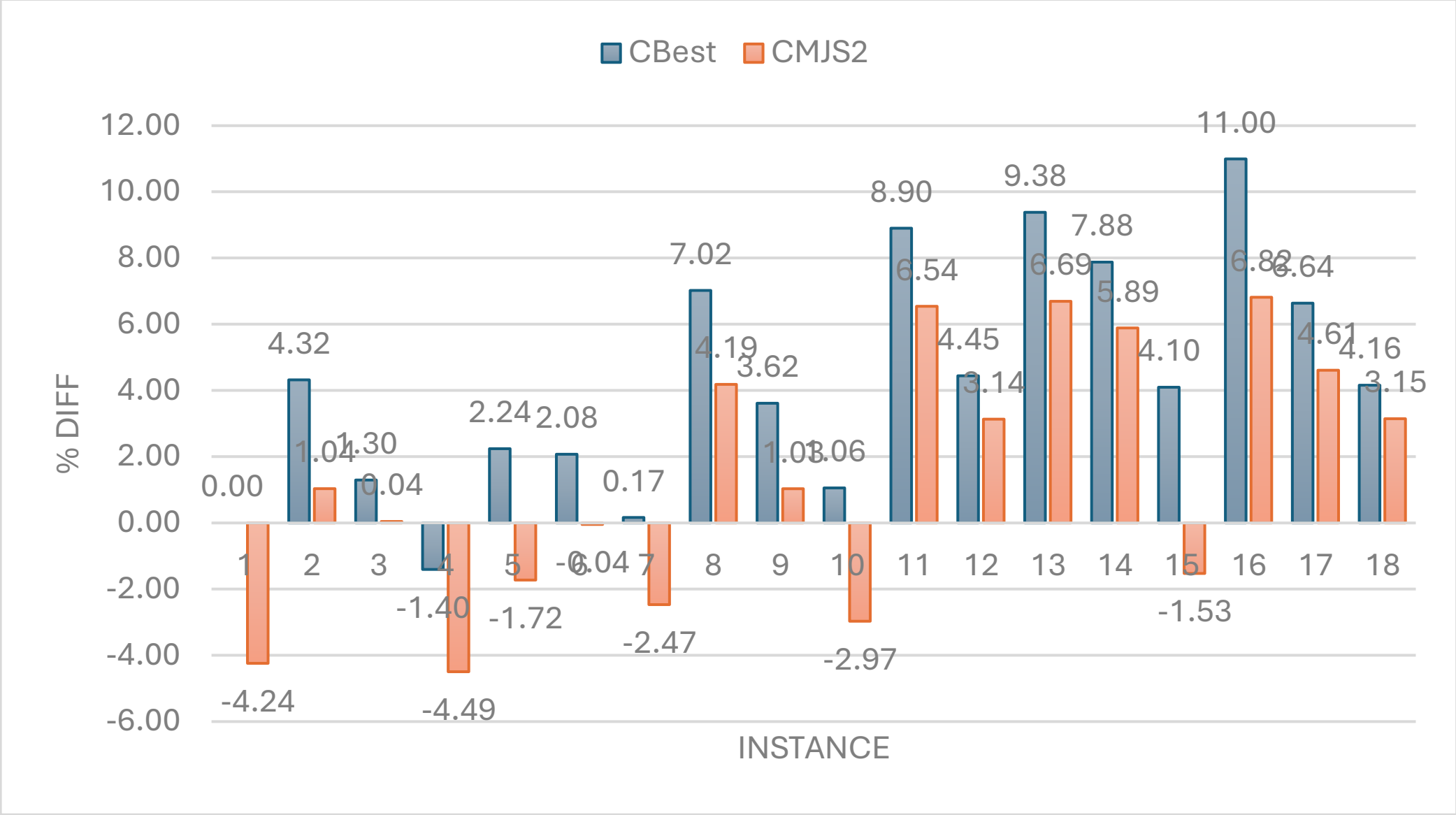

**Figure 4.** Percentage difference of best PISA solution from benchmark solutions (short and long runs)

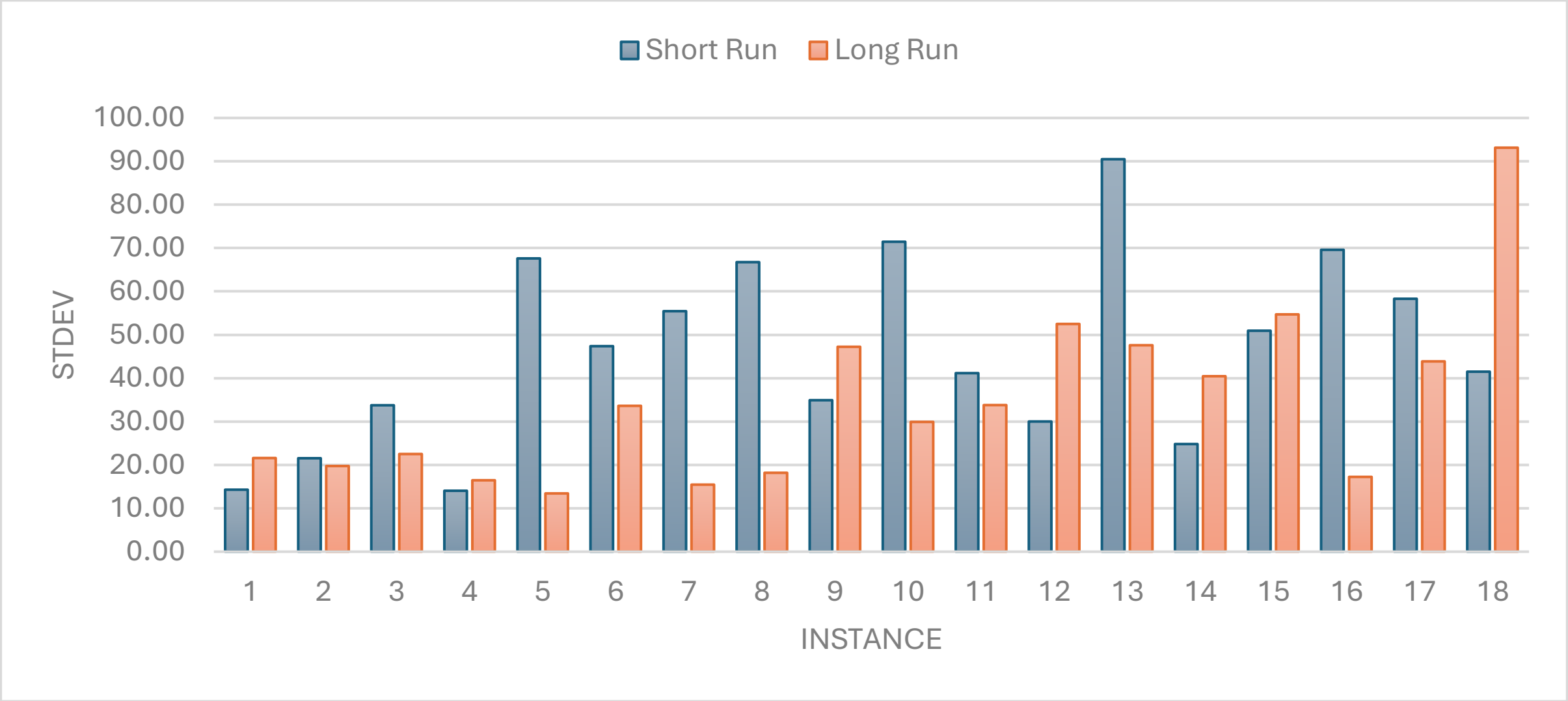

**Figure 5**. Consistency of solutions produced by PISA

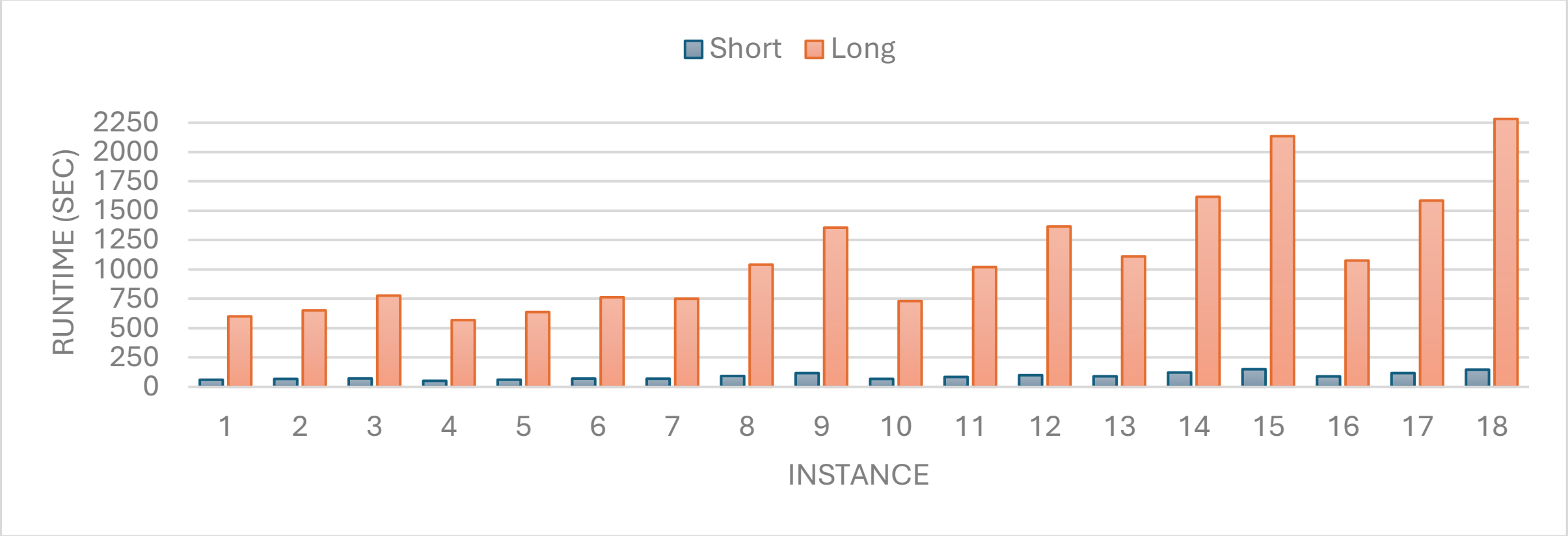

**Figure 6.** Observed differences in PISA runtime

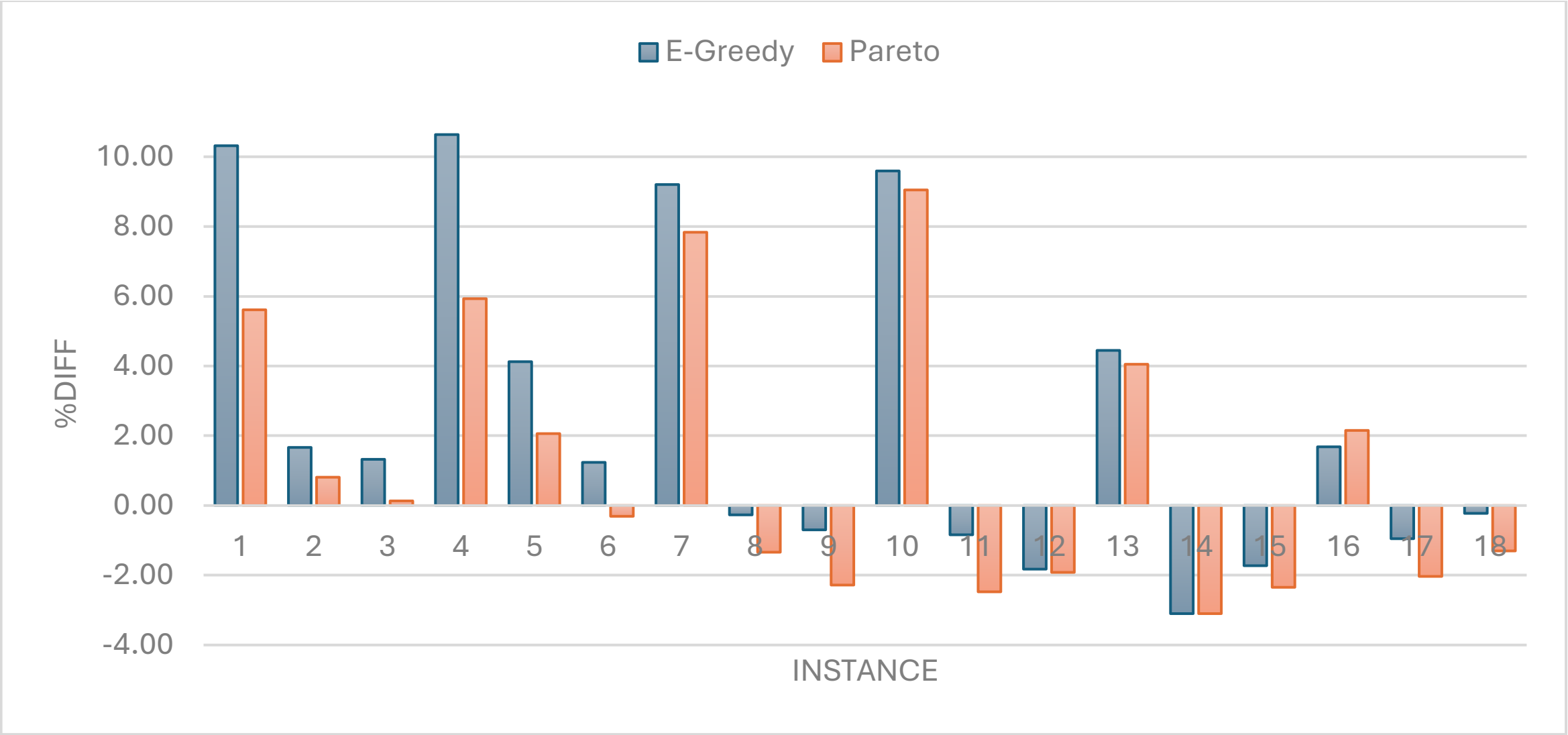

**Figure 7.** QL results - Percentage difference from the best PISA solution

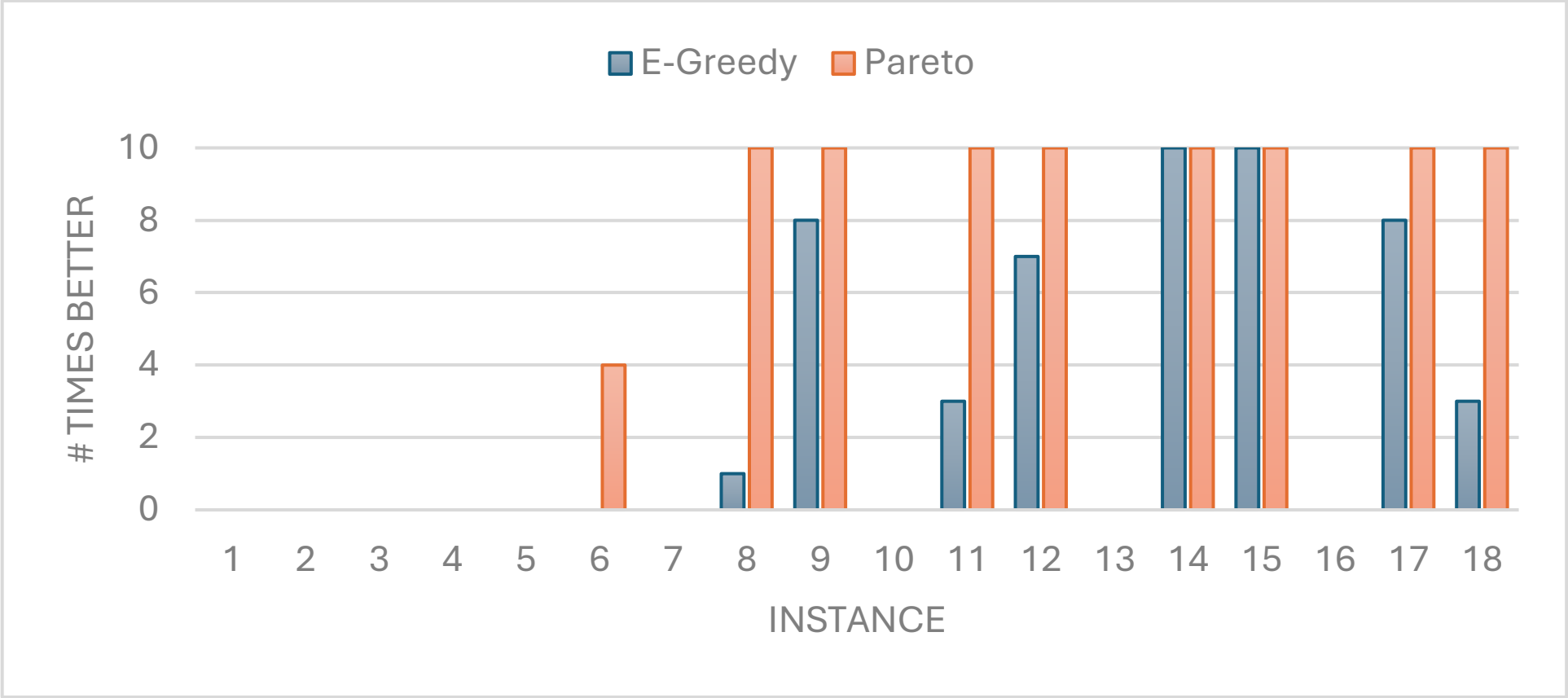

**Figure 8.** Number of trials where QL outperformed PISA

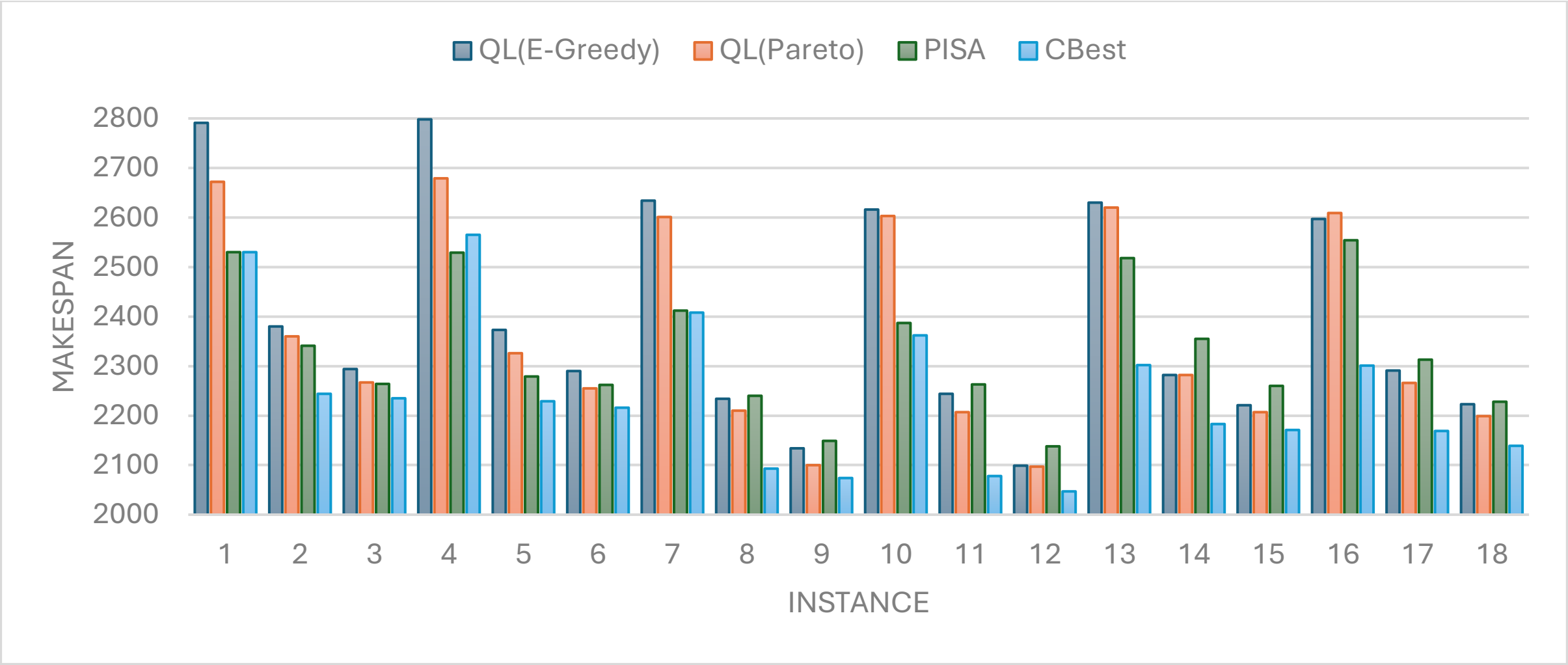

**Figure 9.** Comparison of best makespan

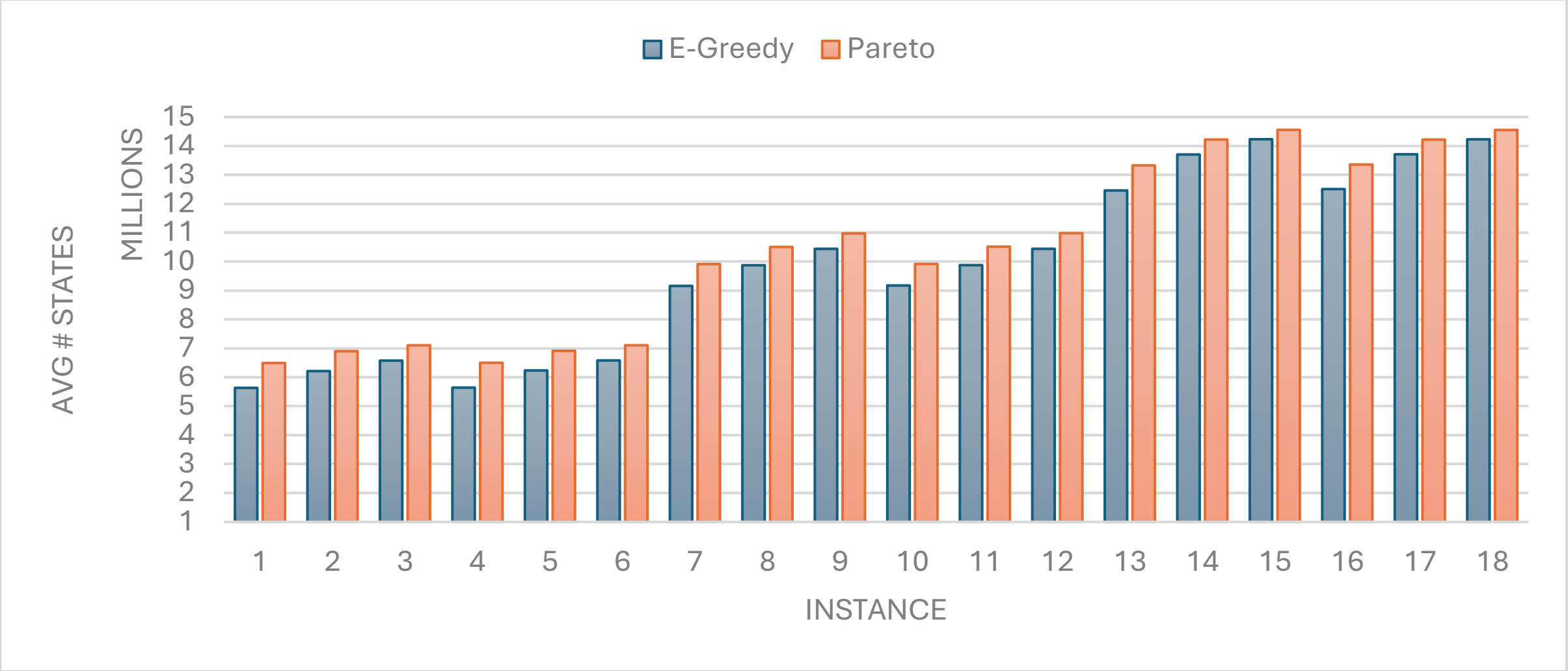

**Figure 10.** Average number of states recorded using QL

As noted earlier, the performance of Q-learning is highly sensitive to the number of training episodes; however, increasing this number also leads to greater memory demands. Although 20,000 episodes proved effective in our experiments, scaling beyond this level becomes increasingly challenging without more sophisticated strategies or access to high-performance computing resources.

To examine whether the previously observed trends would persist, we executed Q-learning on a HPC system using a substantially larger training budget (e.g., 40,000 episodes). The resulting relative improvements are summarised in Table 6. With only a few exceptions, QL at 40,000 episodes consistently achieved lower objective values and reduced variance compared with 20,000 episodes, indicating both improved and more stable performance. The worst-case performance improved markedly.

However, best-case gains were marginal—typically no more than 1%—suggesting that the method may be approaching its effective learning capacity and beginning to exhibit diminishing returns. For illustration, Figure 11 presents the solution quality for instance 18 across varying numbers of episodes. The average number of states, also reported in Table 6, provides an estimate of the associated memory footprint. Relative to the 20,000-episode setting, the memory requirement approximately doubled.

**Table 6.** Relative percentage improvement achieved for QL with 40,000 learning episodes.

| | QL (40000) | | | | | Relative to QL (20000) | | | | Relative to PISA | | | |
|---|---|---|---|---|---|---|---|---|---|---|---|---|---|
| # | min | avg | max | stdev | avg #states | min | avg | max | stdev | min | avg | max | stdev |
| **1** | 2694 | 2719.4 | 2743 | 13.06 | 12,836,274 | 0.82 | -0.44 | -3.45 | -29.38 | 6.48 | 6.35 | 6.40 | -8.63 |
| **2** | 2349 | 2367.8 | 2386 | 10.96 | 13,690,130.9 | **-0.47** | -0.25 | -1.73 | -27.97 | 0.34 | -0.95 | -1.45 | -49.11 |
| **3** | 2263 | 2272.1 | 2281 | 5.384 | 14,123,892.4 | **-0.18** | -0.24 | -1.81 | -47.58 | **-0.04** | -2.32 | -3.63 | -84.05 |
| **4** | 2673 | 2719.4 | 2749 | 25.77 | 12,852,311.5 | **-0.22** | 0.11 | -3.75 | 23.08 | 5.69 | 7.11 | 7.01 | 83.53 |
| **5** | 2321 | 2351.7 | 2372 | 14.8 | 13,717,020.7 | **-0.21** | -0.47 | -1.70 | 41.00 | 1.84 | -1.25 | -4.55 | -78.12 |
| **6** | 2242 | 2257.5 | 2264 | 8.114 | 14,129,897.9 | **-0.58** | -0.35 | -2.46 | -6.10 | **-0.88** | -3.93 | -6.64 | -82.88 |
| **7** | 2608 | 2627.9 | 2640 | 9.848 | 19,725,002.8 | 0.27 | -0.25 | -1.01 | -1.96 | 8.13 | 5.68 | 1.58 | -82.24 |
| **8** | 2192 | 2216.3 | 2228 | 10 | 20,902,809.8 | **-0.81** | -0.31 | -2.79 | -52.16 | **-2.14** | -4.44 | -8.84 | -85.02 |
| **9** | 2095 | 2108.4 | 2116 | 7.058 | 21,915,913.3 | **-0.24** | -0.31 | -2.17 | -16.22 | **-2.51** | -4.20 | -6.12 | -79.79 |
| **10** | 2575 | 2619.3 | 2639 | 18.71 | 19,730,479.4 | **-1.08** | -0.43 | -1.57 | -6.71 | 7.88 | 4.49 | 1.07 | -73.81 |
| **11** | 2211 | 2223.2 | 2231 | 7.51 | 20,926,141.9 | 0.18 | -0.24 | -2.53 | -48.87 | **-2.30** | -5.19 | -7.50 | -81.75 |
| **12** | 2088 | 2098.8 | 2113 | 8.011 | 21,931,316.8 | **-0.43** | -0.39 | -1.63 | -51.99 | **-2.34** | -3.92 | -5.75 | -73.30 |
| **13** | 2598 | 2616.6 | 2640 | 15.28 | 26,471,807.8 | **-0.84** | -0.80 | -2.40 | -40.48 | 3.18 | -2.47 | -5.88 | -83.12 |
| **14** | 2268 | 2282 | 2295 | 8.919 | 28,411,120.8 | **-0.61** | -0.41 | -1.46 | -37.01 | **-3.69** | -4.29 | -5.13 | -64.05 |
| **15** | 2195 | 2209.8 | 2224 | 8.337 | 29,091,023.8 | **-0.54** | -0.28 | -1.33 | -16.54 | **-2.88** | -4.54 | -7.29 | -83.63 |
| **16** | 2583 | 2616.6 | 2643 | 18.77 | 26,544,122.1 | **-1.00** | -0.41 | -1.49 | -32.64 | 1.14 | -2.19 | -4.76 | -73.03 |
| **17** | 2258 | 2274 | 2282 | 6.65 | 28,406,323.4 | **-0.35** | -0.24 | -1.68 | -22.72 | **-2.38** | -4.24 | -8.97 | -88.59 |
| **18** | 2194 | 2208.1 | 2214 | 6.19 | 29,084,958.5 | **-0.23** | -0.07 | -1.25 | -16.11 | **-2.66** | -3.41 | -3.86 | -85.08 |

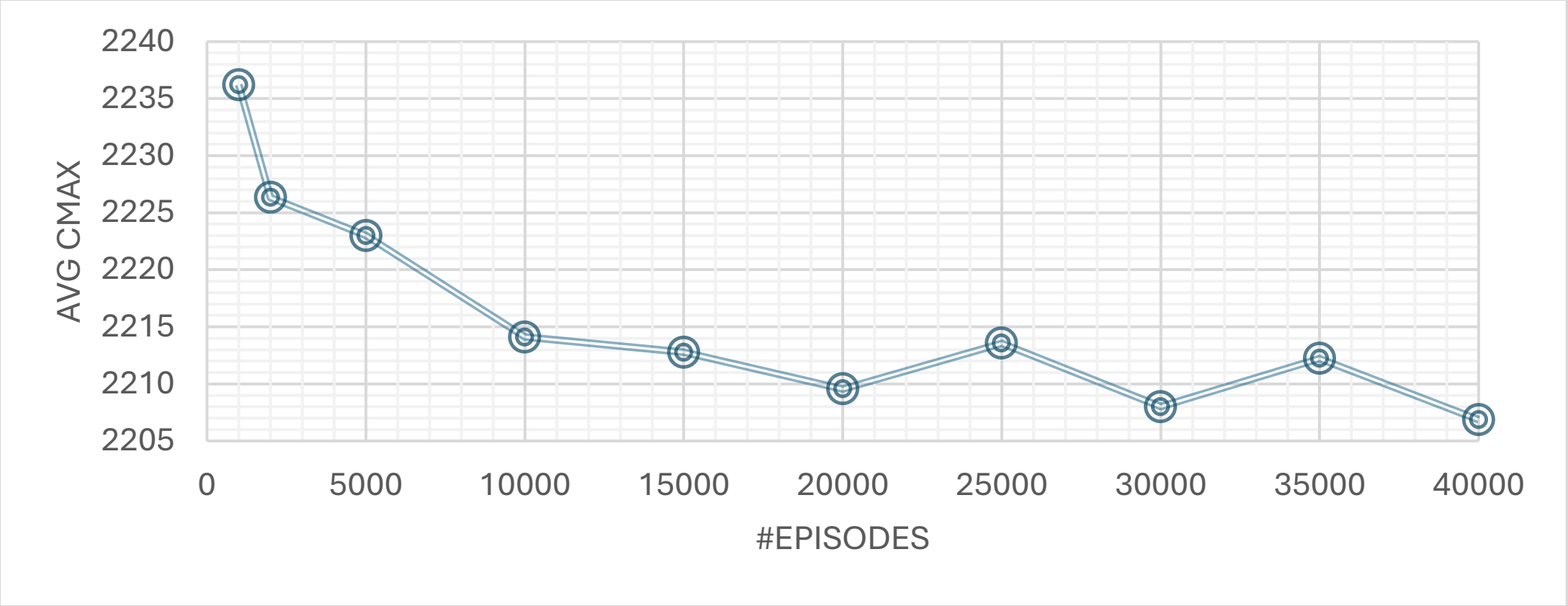

**Figure 11.** QL solution quality relative to the number of episodes (instance 18)

### 4.2. FJSP with Stochastic Processing Times

To further evaluate the two algorithms, the Dauzère-Pérès and Paulli FJSP instances were extended to incorporate stochastic processing times. Specifically, each processing time was modelled as a normally distributed random variable with mean equal to the original deterministic value. The standard deviation was defined as 5% of the expected value divided by 1.96, ensuring that approximately 95% of sampled processing times fall within the interval $[\mu - r, \mu + r]$ where $r = 1.96 \times stdev$.

Preliminary numerical experiments were conducted to identify suitable parameter settings. These tests revealed substantially increased runtimes relative to the deterministic case. Under the configuration $t_{\mathrm{I}} = 10$, $t_{\mathrm{F}} = 0.001$, $maxSteps = 50$, $maxIter = 2000$, the results for the first five instances are reported in Table 7.

**Table 7.** Numerical results of PISA with no parallel simulation evaluations

| Inst. | Cbest | cmax | cpu (#days) | #sim |
|---|---|---|---|---|
| **1** | 2530 | 2593.03 | 1.51 | 1001 |
| **2** | 2244 | 2409.8 | 1.92 | 851 |
| **3** | 2235 | 2289.34 | 2.53 | 814 |
| **4** | 2565 | 2540.1 | 1.5 | 958 |
| **5** | 2229 | 2290.87 | 1.71 | 823 |

Due to the excessive runtimes reported in Table 7, we adopted a single ISA framework combined with parallelised Monte Carlo simulation for solution evaluation. The aggregated results across ten trials and 1,000 scenarios are presented in Table C1 of Appendix C. The choice of 1,000 simulation runs was guided by the convergence-based evaluator, which indicated that both the expected makespan and its variance stabilised after approximately 800 simulations. Leveraging parallel processors for simulation yielded substantial computational gains, reducing runtime from days to hours.

In Table C1, the reported makespan values correspond to expected values. The %RE column represents the relative difference between the minimum makespan in the stochastic setting and that of the deterministic case. Notably, the expected minimum makespan is lower than the deterministic minimum makespan. This outcome arises from the use of normally distributed processing times, which permit sampled values below the mean. Figure 12 illustrates these improvements through a bar-chart comparison.

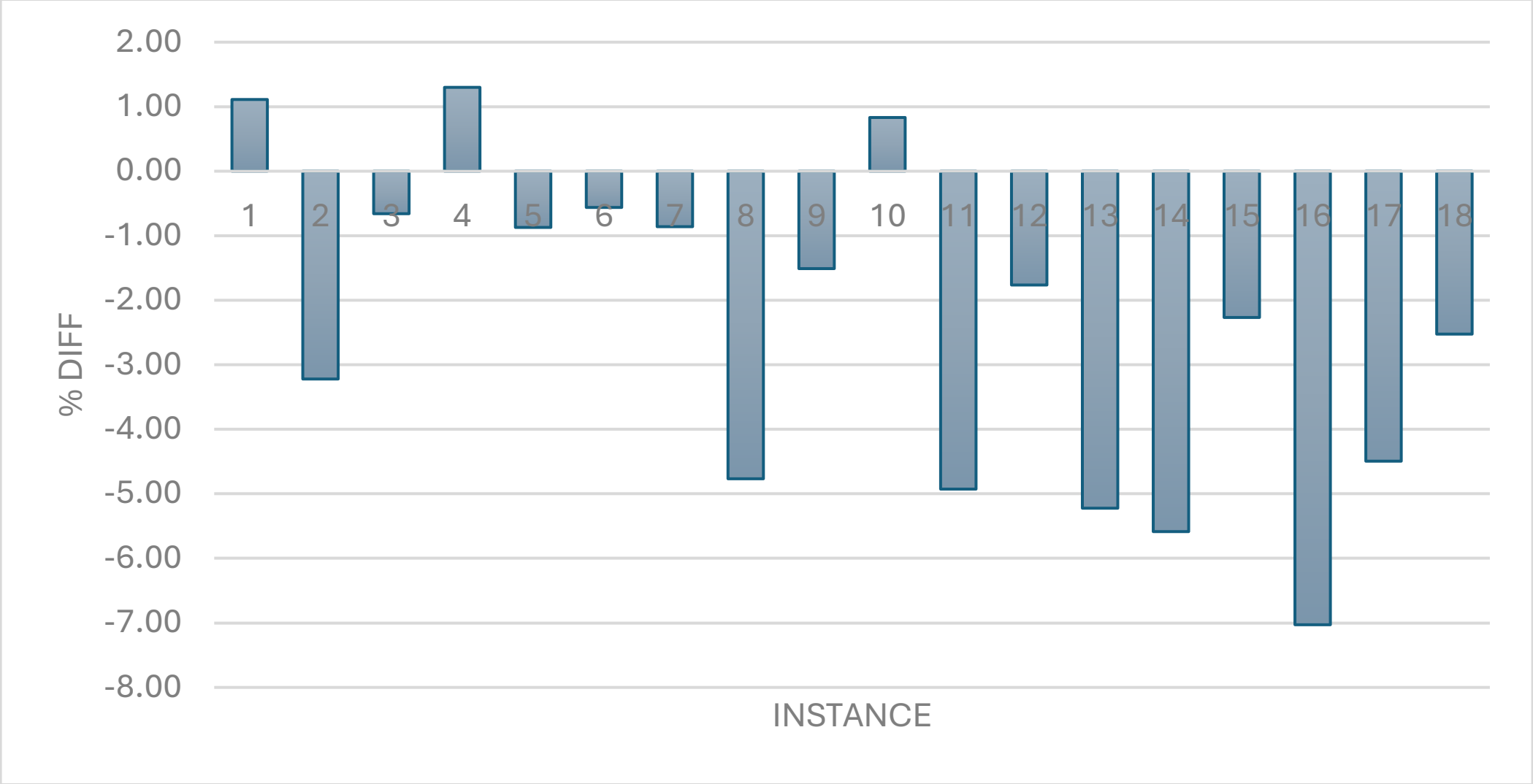

**Figure 12.** Percentage difference of ISA (stochastic) from the best deterministic PISA solution

The Q-learning algorithm was then applied, with results summarised in Appendix C (Tables C2–C5). The %RE columns in Tables C2 and C3 report the relative difference between QL results and the ISA benchmark results presented in Table C1. Pairwise statistical comparisons were conducted to assess performance differences, with results reported in Table 8.

Table 8 indicates that QL (pbs) consistently outperforms QL (ε-greedy) in both best-case and average-case performance. However, ISA achieves superior best-case performance overall, as further illustrated in the bar chart in Figure 13. QL (pbs) becomes competitive primarily in the average-case setting, where its performance is closer to that of ISA. The relatively limited performance of Q-learning can be attributed primarily to the absence of parallel simulation, which restricts the number of episodes that can be evaluated within a given computational budget. As a result, QL is unable to match ISA, which benefits from parallelised simulation and can explore a substantially larger number of candidate solutions simultaneously.

**Table 8.** Pairwise t-tests for stochastic FJSP

| Case | A | B | #better | | mean diff | stdev | two-sided | | | one-sided | |
|---|---|---|---|---|---|---|---|---|---|---|---|
| | | | A | B | B-A | | t | p | sig | p | sig |
| best | QL (egrdy) | QL (pbs) | 3 | 15 | -26.056 | 26.578 | 4.159 | 0.001 | ✓ | 0.000 | ✓ |
| best | ISA | QL (egrdy) | 18 | 0 | 114.601 | 91.366 | 5.322 | 0.000 | ✓ | 0.000 | ✓ |
| best | ISA | QL (pbs) | 18 | 0 | 88.545 | 78.843 | 4.765 | 0.000 | ✓ | 0.000 | ✓ |
| average | QL (egrdy) | QL (pbs) | 0 | 18 | -37.744 | 22.607 | 7.083 | 0.000 | ✓ | 0.000 | ✓ |
| average | ISA | QL (egrdy) | 13 | 5 | 68.963 | 85.270 | 3.431 | 0.003 | ✓ | 0.002 | ✓ |
| average | ISA | QL (pbs) | 9 | 9 | 31.218 | 75.688 | 1.750 | 0.098 | × | 0.049 | ✓ |

Regarding QL, it is noteworthy that the number of observed states exceeds that of the deterministic case; however, it remains within the same order of magnitude. In contrast, the number of terminal states increases substantially, reaching approximately 15,000 compared with around 800 in the deterministic setting. This indicates that stochastic processing times significantly expand the effective state space that QL must explore.

For comparison with the deterministic case, QL was also re-executed using 40,000 learning episodes. The results, summarised in Table 9, show a modest improvement in solution quality relative to ISA. QL now produces superior solutions in four instances, whereas no such improvements were observed previously.

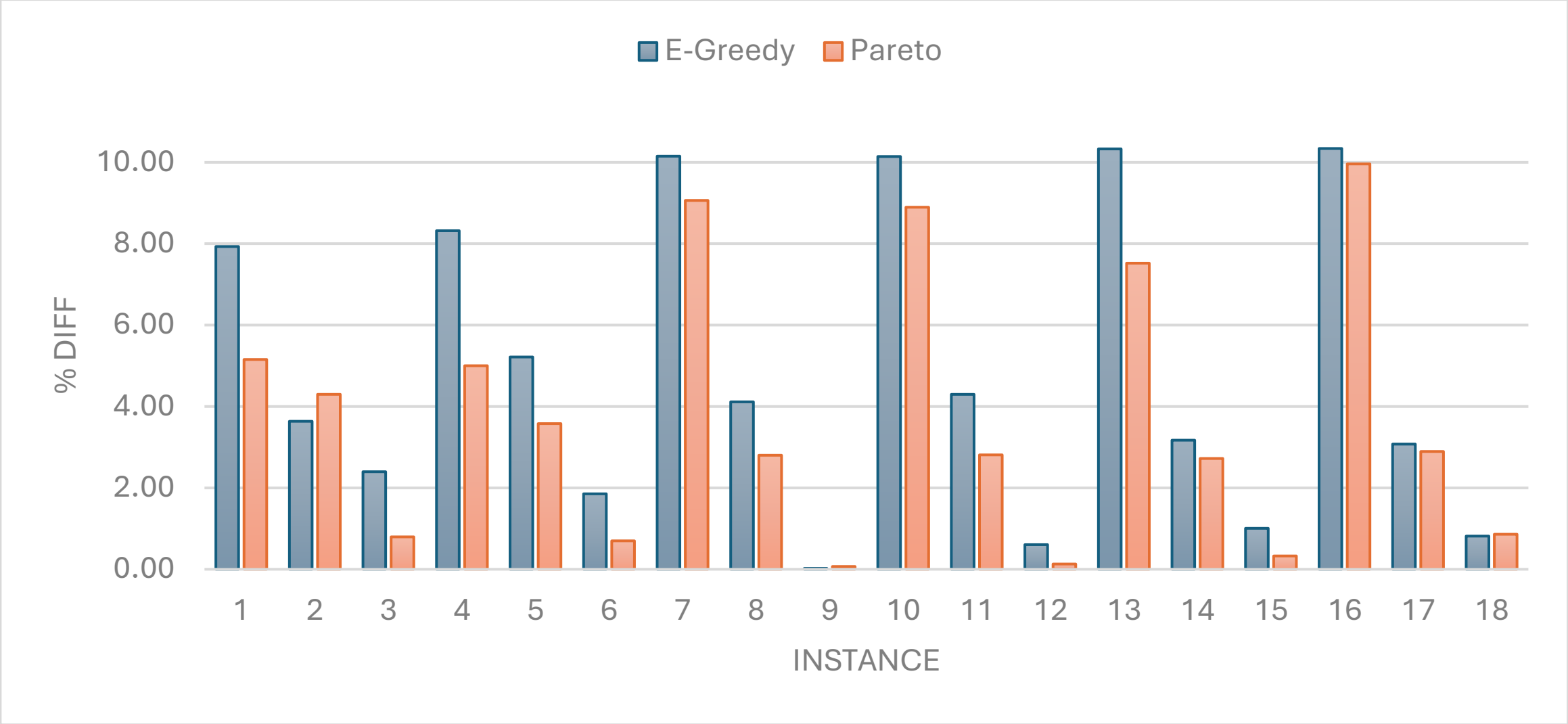

**Figure 13.** Percentage difference from the best ISA solution

**Table 9.** Relative percentage improvement achieved for QL with 40,000 learning episodes.

| | QL (40000) | | | | | Relative to QL (20000) | | | | Relative to ISA | | | |
|---|---|---|---|---|---|---|---|---|---|---|---|---|---|
| **#** | min | avg | max | stdev | avg #states | min | avg | max | stdev | min | avg | max | stdev |
| **1** | 2678 | 2722.2 | 2749 | 20.66 | 14,267,825.2 | **-0.45** | -0.54 | -0.33 | -1.03 | 4.69 | 1.82 | 5.65 | -36.74 |
| **2** | 2353 | 2372.3 | 2381 | 8.01 | 14,638,222.7 | **-0.42** | -0.25 | -0.50 | -24.40 | 3.86 | 2.65 | 3.88 | -46.06 |
| **3** | 2270 | 2276.5 | 2281 | 3.92 | 15,049,760.2 | **0.13** | -0.08 | -0.31 | -45.84 | 0.93 | -1.83 | 0.09 | -79.69 |
| **4** | 2685 | 2714.9 | 2730 | 13.65 | 14,261,773.4 | **-0.19** | -0.51 | -1.02 | -35.06 | 4.81 | 0.15 | 4.78 | -68.45 |
| **5** | 2341 | 2356.1 | 2367 | 8.23 | 14,663,032.3 | 0.04 | -0.11 | -0.25 | -33.42 | 3.62 | 1.55 | 3.21 | -57.73 |
| **6** | 2261 | 2267.8 | 2279 | 6.11 | 15,058,585.7 | **-0.18** | -0.25 | -0.18 | 17.64 | 0.52 | -1.04 | 0.79 | -52.40 |
| **7** | 2595 | 2623.2 | 2641 | 14.36 | 21,013,906 | **-0.50** | -0.29 | -0.49 | -14.50 | 8.52 | 4.91 | 7.69 | -52.95 |
| **8** | 2188 | 2219.8 | 2234 | 14.65 | 21,887,466.4 | **-0.23** | -0.09 | -0.18 | 8.00 | 2.57 | -0.16 | 2.99 | -47.41 |
| **9** | 2111 | 2120 | 2127 | 5.87 | 22,748,647.5 | **-0.33** | -0.26 | -0.56 | -8.76 | **-0.26** | -2.47 | -0.57 | -69.32 |
| **10** | 2601 | 2619.6 | 2648 | 17.44 | 21,030,512.6 | **-0.76** | -0.57 | -0.23 | 63.81 | 8.07 | 6.65 | 8.68 | 4.38 |
| **11** | 2195 | 2224.7 | 2247 | 14.35 | 21,902,330 | **-0.77** | -0.29 | -0.09 | -0.77 | 2.02 | 1.97 | 3.58 | 18.04 |
| **12** | 2095 | 2111.5 | 2124 | 9.12 | 22,757,108.2 | **-0.38** | -0.17 | -0.05 | 31.33 | **-0.25** | -0.63 | 0.54 | -5.21 |
| **13** | 2615 | 2634.6 | 2657 | 12.64 | 28,002,970.3 | 1.91 | 0.21 | -0.15 | -54.99 | 9.58 | 5.75 | 9.62 | -54.73 |
| **14** | 2272 | 2290.9 | 2310 | 10.13 | 29,523,533.3 | **-0.53** | -0.43 | -0.09 | 9.33 | 2.18 | -0.59 | 2.27 | -55.37 |
| **15** | 2203 | 2219.2 | 2229 | 8.31 | 30,305,257.7 | **-0.59** | -0.31 | -0.36 | 31.04 | **-0.26** | -2.59 | -0.18 | -58.99 |
| **16** | 2590 | 2620.1 | 2647 | 20.04 | 28,036,410.1 | **-0.80** | -0.61 | -0.19 | 73.08 | 9.08 | 5.39 | 9.08 | -38.79 |
| **17** | 2266 | 2277.6 | 2283 | 6.10 | 29,540,462.3 | **-0.31** | -0.47 | -0.78 | -29.04 | 2.58 | -2.38 | 0.74 | -83.73 |
| **18** | 2189 | 2213.3 | 2223 | 10.01 | 30,309,953.9 | **-1.22** | -0.33 | -0.22 | 127.63 | **-0.37** | -2.62 | -0.35 | -64.33 |

For test problem 18, a comparative analysis of QL solution quality, in relation to the number of episodes and the deterministic case, is shown in Figure 14.

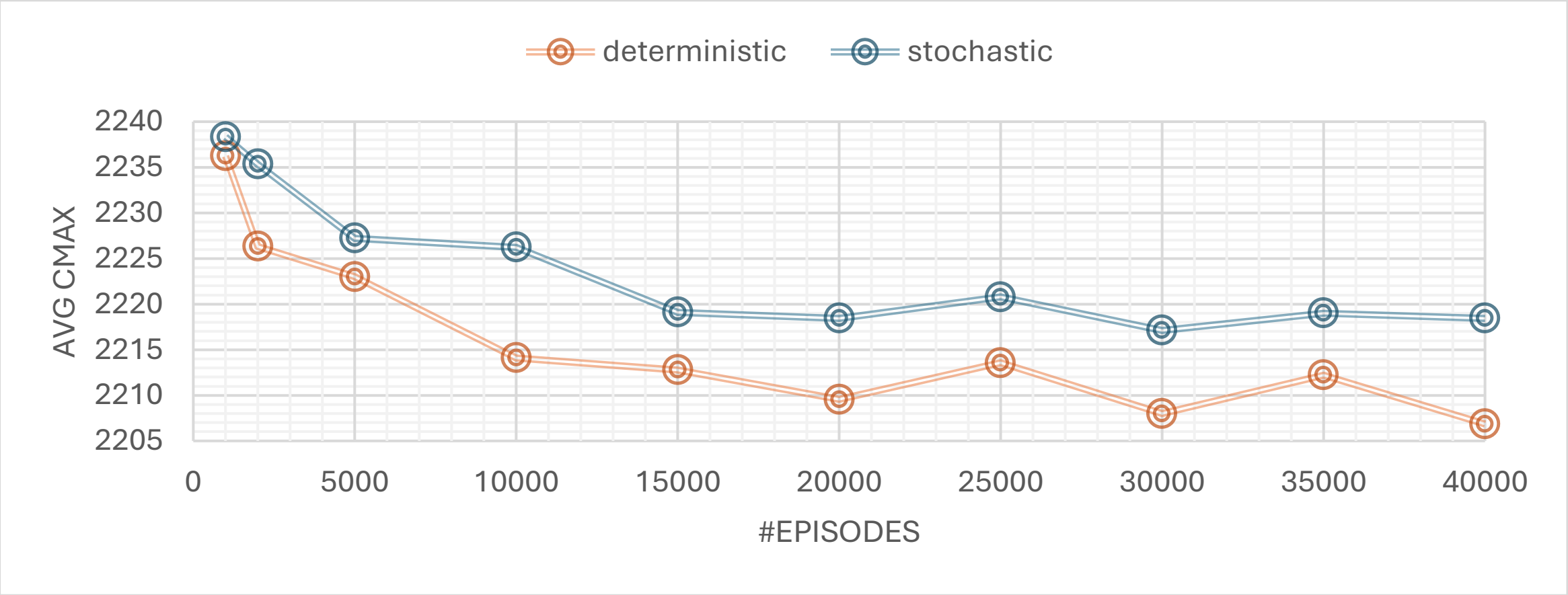

**Figure 14.** Solution quality relative to number of episodes (instance 18).

## 5. Conclusions

In this article, we evaluated key implementation aspects of Simulated Annealing with Intensification (SA) and Q-learning for solving the NP-hard Flexible Job Shop Scheduling Problem. Through extensive numerical experiments, we found that Q-learning delivers competitive performance in terms of both solution quality and runtime, with no substantial degradation as problem size increases.

Overall, the SA variants—enhanced with intensification strategies, parallelisation, and compound perturbation operators—consistently outperformed Q-learning in runtime and generally achieved better best-case solutions, although average performance was broadly comparable. However, attaining this level of SA performance requires considerable domain expertise and careful parameter tuning.

In contrast, Q-learning is conceptually simpler in some respects but more challenging to implement in practice, as it requires a problem-specific encoding of system dynamics for which limited practical guidance is available. In our case, the ISA implementation required less code and was completed in significantly less development time.

The primary limitations of Q-learning in this context arise from state-space explosion, the absence of inherent parallelism, and memory constraints. These issues are not present in ISA, which benefits from parallel evaluation and does not require explicit storage of learned state–action structures. In addition, standard Q-learning is relatively sample inefficient, since each transition is processed sequentially. Collectively, these factors restrict scalability and motivate consideration of more advanced variants, such as asynchronous Q-learning, multi-agent formulations, batch update schemes, and other parallelised learning approaches.

We also investigated several enhancements to the QL framework. Replacing ε-greedy exploration with Pareto-based action selection was straightforward to implement and yielded measurable improvements in FJSP solution quality. However, introducing specialised reward shaping mechanisms did not produce statistically significant gains, suggesting that further effort in this direction is unlikely to be fruitful.

For future research, advanced memory management appears to be a critical area of focus. Pruning the Q-table is a natural direction. Just as human learning involves both retention and forgetting, a controlled forgetting mechanism could be incorporated into Q-learning to mitigate excessive memory growth. The central challenge is how to determine which states should be discarded. One possibility is to base pruning on state utility, for example removing states associated with persistently low-value estimates. Another is to use visitation frequency, eliminating states that are rarely encountered during learning. A hybrid strategy that combines both criteria may also be worth investigating.

Another important direction for future research is batch learning. In a similar spirit to how Simulated Annealing exploits Monte Carlo sampling to evaluate multiple configurations efficiently, analogous batch-based update schemes could be developed for Q-learning to improve both efficiency and scalability.

A related idea is the use of multiple agents within a Q-learning framework. While this is conceptually appealing, implementing a full multi-agent system is non-trivial. Unlike single-agent Q-learning, where the environment is stationary from the learner's perspective, multi-agent settings introduce inherent non-stationarity, as all agents learn and act concurrently. This significantly complicates convergence behaviour, since each agent must account not only for its own actions but also for the evolving policies of others.

A naïve approach based on sharing a single centralised Q-table is therefore unlikely to be effective. It would tend to produce noisy updates due to conflicting experiences, suffer from poor scalability driven by the combinatorial growth of the joint action space, and face fundamental credit assignment difficulties in attributing outcomes to individual agents.

In this article, we focused on the single-objective FJSP. However, real-world applications are typically multi-objective in nature, with competing criteria often playing an equally important role.

Extending the problem to a multi-objective setting introduces substantial challenges for both Q-learning and Simulated Annealing, particularly in terms of increased computational burden and memory requirements. A systematic investigation of which adaptations are effective—and which are not—would therefore be a valuable contribution.

In a multi-objective Q-learning framework, both the reward structure and Q-values would need to be extended, for example through vector-valued representations or by introducing an additional objective index. A key difficulty in this setting lies in managing the growing set of trade-off solutions, particularly in determining when non-dominated actions should be retained or discarded in memory. Developing principled strategies for this selection process remains an open and important challenge.

**Disclosure statement**

No potential conflict of interest was reported by the author(s).

**Funding**
This research was funded by the United States Air Force, Asian Office of Aerospace Research and Development.

## Appendix A: Implementation Details – Class Details

| Class | Attributes | Main Functions |
|---|---|---|
| SystemState (state of the system at a specific time) | ▪ Time t;  // System Time<br>▪ unordered_map<string, StateRecord> records;<br>where:<br>class StateRecord {<br>ObjState state; string more_info; Object *obj;<br>} | ▪ template <typename T> void add_more_info(string& key, T value);<br>▪ template <typename T> T get_more_info(const string& key);<br>▪ void update(string name, StateRecord& record);<br>▪ void add(StateRecord &record);<br>▪ void log(ostream& os);<br>▪ StateRecord* lookup_record(string name);<br>▪ std::any lookup_info(string name); |
| SystemHistory | ▪ vector<SystemState*> states;<br>▪ vector<Action*> actions;<br>▪ System* system; | ▪ log();<br>▪ record_state(SystemState *state);<br>▪ record_action(Action &action); |
| System | ▪ SystemHistory history;<br>▪ SystemState* state_zero;<br>▪ SystemState* current_state; | ▪ virtual SystemState* extract_state(Time t);<br>▪ virtual void get_possible_actions(SystemState& state, Actions& actions) = 0;<br>▪ virtual void take_action(SystemState &state, Action &action) = 0;<br>▪ virtual void transition(SystemState &current, Time t) = 0;<br>▪ virtual double compute_reward(SystemState& previous, SystemState& next, Action action) = 0;<br>▪ virtual Time* get_next_decision_point(SystemState& state, bool exclude_current = false) = 0;<br>▪ virtual string hash_state(SystemState &state) = 0; |
| Action | ▪ std::any info;<br>▪ double reward;<br>▪ Time when; | ▪ add_reward(reward);<br>▪ template<typename T> T get();<br>▪ template<typename T> define(T info); |
| Object | ▪ string descr;<br>▪ ObjState state; | ▪ void set_state(); |
| States | ▪ string | |
| Job: Object | ▪ int num;<br>▪ string descr;<br>▪ JobShop* shop;<br>▪ vector<Operation*> operations;<br>▪ Time rlt;<br>▪ Operation* current; | ▪ void add_operation(Operation& operation);<br>▪ void log_details(ostream& os);<br>▪ void update_state(ObjState::Values value);<br>▪ Operation* first();<br>▪ Operation* last();<br>▪ Operation* succ(Operation* operation);<br>▪ Operation* pred(Operation* operation); |
| Operation: Object | ▪ int id, num;<br>▪ Job& job;<br>▪ Time start, end;<br>▪ Parameter* pt;<br>▪ vector<Machine*> options;<br>▪ unordered_map<string, Parameter*> pts;<br>▪ Machine *machine; | ▪ void add_option(Machine* _mach, Parameter* pt);<br>▪ void log_details(ostream& os);<br>▪ void update_state(ObjState::Values value);<br>▪ Time get_remain(Time t);<br>▪ void assign_machine(Machine* mach);<br>▪ Time get_pt(Machine* machine, bool sample=true);<br>▪ Operation* succ();<br>▪ Operation* pred(); |
| Machine: Object | ▪ string descr;<br>▪ int num;<br>▪ JobShop* shop; | ▪ void update_state(ObjState::Values value); |
| JobShop: System | ▪ using PriorityOrder = vector<Assignment>;<br>▪ vector<Job*> jobs;<br>▪ vector<Machine*> machines;<br>▪ vector<Operation*> operations;<br>▪ vector<Assignment*> feasible_assignments;<br>▪ unordered_map<string, int> assignment_index;<br>▪ Time cmax;<br>▪ Schedule* schedule; | ▪ void define_action_space();<br>▪ void schedule_first_operation(Job *job);<br>▪ double evaluate_actions(vector<Action>& actions) override;<br>▪ Metrics* evaluate_schedule(Schedule& schedule);<br>▪ Metrics* evaluate(PriorityOrder& order, bool stoch=true);<br>▪ Metrics* evaluate_expected_performance(PriorityOrder& order, int nb_scen = 1);<br>▪ PriorityOrder* extract_priority_order();<br>▪ PriorityOrder* get_priority_order(vector<Action>& actions);<br>▪ MachineSequences* extract_sequences(PriorityOrder& order);<br>▪ Schedule* build_schedule(PriorityOrder &order, bool sample); |

## Appendix B

**Table B1.** Numerical results of PISA over 10 "shorter" trials

| | **cmax** | | | | **cpu (sec)** | | | | **%DIFF** | |
|---|---|---|---|---|---|---|---|---|---|---|
| **Inst.** | min | avg | max | stdev | min | avg | max | stdev | Cbest | CMSJ2 |
| **1** | 2583 | 2650 | 2615.5 | 21.62 | 51 | 69 | 60.8 | 5.96 | 2.09 | -2.23 |
| **2** | 2294 | 2363 | 2341.2 | 24.46 | 65 | 72 | 67.9 | 2.23 | 2.23 | -0.99 |
| **3** | 2304 | 2408 | 2370.1 | 30.90 | 66 | 75 | 71.8 | 2.70 | 3.09 | 1.81 |
| **4** | 2564 | 2694 | 2614 | 38.40 | 49 | 54 | 51.6 | 1.90 | -0.04 | -3.17 |
| **5** | 2245 | 2367 | 2324.2 | 33.38 | 57 | 64 | 61.4 | 2.27 | 0.72 | -3.19 |
| **6** | 2276 | 2401 | 2340.2 | 40.30 | 65 | 74 | 70.6 | 2.88 | 2.71 | 0.57 |
| **7** | 2434 | 2526 | 2472.7 | 31.02 | 65 | 73 | 69.7 | 2.71 | 1.08 | -1.58 |
| **8** | 2291 | 2439 | 2391.2 | 45.92 | 89 | 96 | 92.7 | 2.45 | 9.46 | 6.56 |
| **9** | 2187 | 2347 | 2229.1 | 45.81 | 113 | 126 | 117.3 | 3.53 | 5.45 | 2.82 |
| **10** | 2433 | 2498 | 2465.6 | 22.14 | 64 | 71 | 68.2 | 1.99 | 3.01 | -1.10 |
| **11** | 2268 | 2444 | 2370.7 | 56.11 | 78 | 100 | 84.7 | 6.25 | 9.14 | 6.78 |
| **12** | 2149 | 2305 | 2224.5 | 47.72 | 92 | 107 | 99.8 | 5.45 | 4.98 | 3.67 |
| **13** | 2433 | 2539 | 2491.4 | 32.89 | 85 | 94 | 90.1 | 2.77 | 5.69 | 3.09 |
| **14** | 2414 | 2481 | 2440.5 | 23.48 | 118 | 129 | 122.5 | 3.37 | 10.58 | 8.54 |
| **15** | 2230 | 2396 | 2316.7 | 51.04 | 145 | 155 | 150.4 | 3.20 | 2.72 | -2.83 |
| **16** | 2432 | 2513 | 2470.4 | 24.91 | 87 | 91 | 89.3 | 1.34 | 5.69 | 1.71 |
| **17** | 2321 | 2545 | 2428.3 | 79.62 | 88 | 122 | 117.5 | 10.41 | 7.01 | 4.98 |
| **18** | 2273 | 2408 | 2331.3 | 44.16 | 145 | 151 | 146.7 | 1.89 | 6.26 | 5.23 |
| | | | | | | | | Avg. | 4.55 | 1.70 |

**Table B2.** Numerical results of PISA over 10 "longer" trials

| | **cmax** | | | | **cpu (sec)** | | | | **%DIFF** | |
|---|---|---|---|---|---|---|---|---|---|---|
| **Inst.** | min | avg | max | stdev | min | avg | max | stdev | Cbest | CMSJ2 |
| **1** | 2530 | 2578 | 2557 | 14.29 | 551.86 | 621.21 | 600.47 | 21.59 | 0.00 | -4.24 |
| **2** | 2341 | 2421 | 2390.6 | 21.54 | 611.22 | 675.99 | 651.65 | 19.75 | 4.32 | 1.04 |
| **3** | 2264 | 2367 | 2326 | 33.76 | 732.53 | 808.89 | 777.62 | 22.51 | 1.30 | 0.04 |
| **4** | 2529 | 2569 | 2538.9 | 14.04 | 537.58 | 588.59 | 568.80 | 16.46 | -1.40 | -4.49 |
| **5** | 2279 | 2485 | 2381.5 | 67.61 | 623.26 | 660.99 | 637.46 | 13.43 | 2.24 | -1.72 |
| **6** | 2262 | 2425 | 2349.8 | 47.38 | 714.21 | 816.18 | 763.35 | 33.62 | 2.08 | -0.04 |
| **7** | 2412 | 2599 | 2486.7 | 55.44 | 727.40 | 778.40 | 750.88 | 15.45 | 0.17 | -2.47 |
| **8** | 2240 | 2444 | 2319.2 | 66.75 | 1009.81 | 1065.61 | 1041.23 | 18.18 | 7.02 | 4.19 |
| **9** | 2149 | 2254 | 2200.9 | 34.92 | 1268.29 | 1415.44 | 1356.24 | 47.24 | 3.62 | 1.03 |
| **10** | 2387 | 2611 | 2506.7 | 71.45 | 700.75 | 809.16 | 730.97 | 29.92 | 1.06 | -2.97 |
| **11** | 2263 | 2412 | 2345 | 41.16 | 966.65 | 1057.28 | 1020.29 | 33.80 | 8.90 | 6.54 |
| **12** | 2138 | 2242 | 2184.4 | 30.01 | 1280.69 | 1444.48 | 1366.35 | 52.50 | 4.45 | 3.14 |
| **13** | 2518 | 2805 | 2683 | 90.49 | 1044.45 | 1219.60 | 1111.26 | 47.59 | 9.38 | 6.69 |
| **14** | 2355 | 2419 | 2384.2 | 24.81 | 1556.79 | 1665.71 | 1619.14 | 40.45 | 7.88 | 5.89 |
| **15** | 2260 | 2399 | 2314.9 | 50.94 | 2055.73 | 2188.35 | 2135.41 | 54.71 | 4.10 | -1.53 |
| **16** | 2554 | 2775 | 2675.2 | 69.60 | 1048.26 | 1101.07 | 1076.27 | 17.24 | 11.00 | 6.82 |
| **17** | 2313 | 2507 | 2374.6 | 58.31 | 1533.35 | 1641.76 | 1587.12 | 43.87 | 6.64 | 4.61 |
| **18** | 2228 | 2369 | 2286 | 41.50 | 2136.76 | 2234.01 | 2282.50 | 93.14 | 4.16 | 3.15 |
| | | | | | | | | Avg. | 4.27 | 1.43 |

**Table B3.** Summarised numerical results of QL (E-greedy exploration)

| | **cmax** | | | | | **cpu (s)** | | | | | |
|---|---|---|---|---|---|---|---|---|---|---|---|
| **Inst.** | min | avg | max | stdev | %PISA | min | avg | max | stdev | %PISA | Ratio |
| **1** | 2791 | 2816 | 2841 | 18.49 | **10.32** | 1097 | 1110 | 1143 | 12.84 | 84.77 | **1.85** |
| **2** | 2380 | 2406 | 2428 | 15.22 | **1.67** | 1183 | 1198 | 1208 | 8.58 | 83.78 | **1.84** |
| **3** | 2294 | 2310 | 2323 | 10.27 | **1.33** | 1269 | 1284 | 1303 | 11.49 | 65.17 | **1.65** |
| **4** | 2798 | 2830 | 2856 | 20.94 | **10.64** | 1098 | 1110 | 1119 | 7.88 | 95.06 | **1.95** |
| **5** | 2373 | 2390 | 2413 | 10.49 | **4.12** | 1175 | 1189 | 1204 | 10.79 | 86.59 | **1.87** |
| **6** | 2290 | 2301 | 2321 | 8.64 | **1.24** | 1253 | 1275 | 1290 | 10.20 | 66.97 | **1.67** |
| **7** | 2634 | 2654 | 2667 | 10.04 | **9.20** | 2644 | 2691 | 2779 | 39.20 | 258.35 | **3.58** |
| **8** | 2234 | 2264 | 2292 | 20.90 | **-0.27** | 2922 | 2973 | 3018 | 35.14 | 185.56 | **2.86** |
| **9** | 2134 | 2144 | 2163 | 8.43 | **-0.70** | 3349 | 3392 | 3464 | 36.92 | 150.12 | **2.5** |
| **10** | 2616 | 2657 | 2681 | 20.06 | **9.59** | 2643 | 2687 | 2736 | 29.79 | 267.59 | **3.68** |
| **11** | 2244 | 2269 | 2289 | 14.69 | **-0.84** | 2922 | 2948 | 2991 | 20.34 | 188.94 | **2.89** |
| **12** | 2099 | 2125 | 2148 | 16.69 | **-1.82** | 3334 | 3379 | 3446 | 36.62 | 147.33 | **2.47** |
| **13** | 2630 | 2673 | 2705 | 25.66 | **4.45** | 4763 | 4808 | 4938 | 53.80 | 332.68 | **4.33** |
| **14** | 2282 | 2314 | 2329 | 14.16 | **-3.10** | 5530 | 5593 | 5690 | 47.14 | 245.44 | **3.45** |
| **15** | 2221 | 2238 | 2254 | 9.99 | **-1.73** | 6548 | 6647 | 6760 | 65.14 | 211.25 | **3.11** |
| **16** | 2597 | 2660 | 2683 | 27.86 | **1.68** | 4808 | 4863 | 4926 | 40.91 | 351.82 | **4.52** |
| **17** | 2291 | 2305 | 2321 | 8.60 | **-0.95** | 5485 | 5553 | 5634 | 43.70 | 249.90 | **3.5** |
| **18** | 2223 | 2234 | 2242 | 7.38 | **-0.22** | 6596 | 6648 | 6688 | 34.60 | 191.26 | **2.91** |

**Table B4.** Summarised numerical results of QL (Pareto based selection)

| | **cmax** | | | | | | **cpu (s)** | | | | | |
|---|---|---|---|---|---|---|---|---|---|---|---|---|
| **Inst.** | min | avg | max | stdev | %PISA | %E-GRY | min | avg | max | stdev | %PISA | Ratio |
| **1** | 2672 | 2732 | 2749 | 22.59 | **5.61** | -3.01 | 1547 | 1580 | 1605 | 17.51 | 163.1 | **2.63** |
| **2** | 2360 | 2374 | 2400 | 14.09 | **0.81** | -1.33 | 1794 | 1831 | 1921 | 47.20 | 181 | **2.81** |
| **3** | 2267 | 2278 | 2285 | 5.583 | **0.13** | -1.42 | 2006 | 2028 | 2045 | 9.99 | 160.7 | **2.61** |
| **4** | 2679 | 2716 | 2742 | 23.5 | **5.93** | -4.00 | 1565 | 1580 | 1594 | 10.37 | 177.8 | **2.78** |
| **5** | 2326 | 2363 | 2383 | 16.28 | **2.06** | -1.14 | 1780 | 1809 | 1847 | 25.22 | 183.7 | **2.84** |
| **6** | 2255 | 2265 | 2279 | 7.905 | **-0.31** | -1.55 | 2006 | 2032 | 2065 | 23.71 | 166.2 | **2.66** |
| **7** | 2601 | 2634 | 2648 | 13.87 | **7.84** | -0.73 | 3749 | 3795 | 3856 | 28.14 | 405.4 | **5.05** |
| **8** | 2210 | 2223 | 2235 | 6.935 | **-1.34** | -1.80 | 4819 | 4861 | 4935 | 34.84 | 366.9 | **4.67** |
| **9** | 2100 | 2115 | 2125 | 7.76 | **-2.28** | -1.35 | 8281 | 8409 | 8620 | 95.54 | 520 | **6.2** |
| **10** | 2603 | 2631 | 2653 | 14.32 | **9.05** | -0.98 | 3742 | 3783 | 3847 | 33.47 | 417.5 | **5.18** |
| **11** | 2207 | 2229 | 2249 | 13.94 | **-2.47** | -1.77 | 4735 | 4795 | 4860 | 45.05 | 370 | **4.7** |
| **12** | 2097 | 2107 | 2116 | 7.218 | **-1.92** | -0.83 | 8266 | 8408 | 8708 | 130.37 | 515.4 | **6.15** |
| **13** | 2620 | 2638 | 2668 | 16.17 | **4.05** | -1.31 | 7434 | 7524 | 7606 | 55.08 | 577.1 | **6.77** |
| **14** | 2282 | 2291 | 2302 | 6.931 | **-3.10** | -0.96 | 13071 | 13152 | 13316 | 79.87 | 712.3 | **8.12** |
| **15** | 2207 | 2216 | 2223 | 4.807 | **-2.35** | -0.97 | 27104 | 27377 | 27556 | 129.89 | 1182 | **12.8** |
| **16** | 2609 | 2628 | 2651 | 12.86 | **2.15** | -1.24 | 7481 | 7516 | 7568 | 25.39 | 598.3 | **6.98** |
| **17** | 2266 | 2280 | 2290 | 7.412 | **-2.03** | -1.12 | 12908 | 13107 | 13453 | 152.97 | 725.8 | **8.26** |
| **18** | 2199 | 2210 | 2218 | 5.87 | **-1.3** | -1.07 | 27121 | 27514 | 27887 | 285.32 | 1105 | **12.1** |

**Table B5.** Number of states recorded and number of full solutions (E-greedy)

| | **#States** | | | | **#Solutions** | | | |
|---|---|---|---|---|---|---|---|---|
| **Inst.** | min | avg | max | stdev | min | avg | max | stdev |
| **1** | 5,620,791 | 5,633,970.2 | 5,641,356 | 5935.65 | 829 | 846.4 | 865 | 11.13 |
| **2** | 6,188,410 | 6,213,412.9 | 6,224,461 | 11140.23 | 819 | 832.5 | 842 | 6.69 |
| **3** | 6,567,142 | 6,575,572 | 6,582,309 | 4892.69 | 811 | 819.9 | 827 | 5.04 |
| **4** | 5,635,452 | 5,642,794.5 | 5,655,284 | 5986.03 | 809 | 829.8 | 847 | 12.91 |
| **5** | 6,218,937 | 6,236,259 | 6,244,719 | 8428.29 | 813 | 825.6 | 833 | 6.88 |
| **6** | 6,578,664 | 6,582,489.3 | 6,589,601 | 3345.38 | 803 | 812.3 | 824 | 6.18 |
| **7** | 9,152,056 | 9,158,830.6 | 9,164,070 | 3632.65 | 1128 | 1147.7 | 1164 | 11.14 |
| **8** | 9,858,582 | 9,874,313.1 | 9,890,482 | 11147.91 | 750 | 774.3 | 805 | 17.75 |
| **9** | 10,435,151 | 10,441,830.1 | 10,451,593 | 6054.76 | 655 | 664.7 | 676 | 5.64 |
| **10** | 9,164,621 | 9,173,829.5 | 9,191,757 | 7865.47 | 1122 | 1128.5 | 1137 | 4.60 |
| **11** | 9,861,388 | 9,877,241.2 | 9,898,679 | 12885.58 | 750 | 775 | 797 | 18.00 |
| **12** | 10,435,748 | 10,441,476.8 | 10,450,683 | 5742.95 | 651 | 660.8 | 677 | 8.32 |
| **13** | 12,440,432 | 12,459,832.1 | 12,473,007 | 10867.32 | 937 | 946.7 | 962 | 7.45 |
| **14** | 13,690,476 | 13,703,308.9 | 13,723,950 | 9974.26 | 600 | 612.9 | 622 | 7.56 |
| **15** | 14,218,159 | 14,229,919.3 | 14,238,755 | 7356.11 | 553 | 563.9 | 579 | 9.45 |
| **16** | 12,491,109 | 12,508,239.7 | 12,524,370 | 10579.50 | 900 | 912.5 | 932 | 11.03 |
| **17** | 13,695,768 | 13,710,162.1 | 13,725,147 | 9518.68 | 595 | 605.9 | 621 | 8.18 |
| **18** | 14,219,585 | 14,228,323.8 | 14,243,054 | 7191.30 | 532 | 538.4 | 550 | 5.42 |

**Table B6.** Number of states recorded and number of full solutions (Pareto based selection)

| | **#States** | | | | **#Solutions** | | | |
|---|---|---|---|---|---|---|---|---|
| **Inst.** | min | avg | max | stdev | min | avg | max | stdev |
| **1** | 6,492,726 | 6,494,262.4 | 6,495,486 | 757.14 | 708 | 716.6 | 728 | 7.31 |
| **2** | 6,901,516 | 6,902,729.3 | 6,903,935 | 995.17 | 700 | 714.6 | 725 | 6.45 |
| **3** | 7,103,975 | 7,105,687.9 | 7,107,351 | 1117.99 | 632 | 645.6 | 661 | 9.89 |
| **4** | 6,499,494 | 6,500,766.3 | 6,502,373 | 1031.79 | 713 | 725.5 | 734 | 6.35 |
| **5** | 6,913,843 | 6,915,429.8 | 6,916,247 | 680.86 | 695 | 707.3 | 718 | 7.21 |
| **6** | 7,106,408 | 7,107,434.3 | 7,108,775 | 800.97 | 651 | 663.1 | 673 | 7.62 |
| **7** | 9,911,708 | 9,913,520 | 9,914,556 | 894.26 | 771 | 778.7 | 785 | 5.64 |
| **8** | 10,504,151 | 10,506,552.1 | 10,508,124 | 1153.64 | 592 | 601.3 | 618 | 9.27 |
| **9** | 10,974,337 | 10,975,401.5 | 10,976,701 | 896.80 | 571 | 581 | 592 | 5.48 |
| **10** | 9,915,997 | 9,916,602.7 | 9,917,700 | 495.69 | 773 | 781.7 | 798 | 7.89 |
| **11** | 10,515,307 | 10,516,812 | 10,518,608 | 1093.16 | 596 | 608.2 | 617 | 6.96 |
| **12** | 10,982,244 | 10,982,764.1 | 10,983,468 | 473.23 | 568 | 575.4 | 582 | 4.25 |
| **13** | 13,322,949 | 13,325,199.5 | 13,328,157 | 1588.61 | 869 | 883 | 900 | 8.12 |
| **14** | 14,219,231 | 14,220,131.3 | 14,221,807 | 886.54 | 607 | 648.9 | 671 | 21.38 |
| **15** | 14,553,757 | 14,555,819.8 | 14,556,916 | 946.55 | 604 | 666.1 | 868 | 106.41 |
| **16** | 13,354,910 | 13,356,276.5 | 13,358,620 | 1196.57 | 646 | 829 | 888 | 94.66 |
| **17** | 14,216,298 | 14,217,894 | 14,221,392 | 1643.16 | 620 | 642 | 661 | 12.68 |
| **18** | 14,551,321 | 14,552,310.8 | 14,553,268 | 649.33 | 616 | 622.5 | 636 | 5.89 |

## Appendix C

**Table C1.** Numerical results of single ISA with parallel simulation evaluations

| | cmax | | | | | cpu (hrs) | | | |
|---|---|---|---|---|---|---|---|---|---|
| **Inst.** | min | avg | max | stdev | %RE DET | min | max | avg | std |
| **1** | 2558.11 | 2673.50 | 2602.06 | 32.66 | **1.11** | 4.58 | 7.50 | 5.43 | 1.05 |
| **2** | 2265.61 | 2311.03 | 2292.01 | 14.86 | **-3.22** | 4.18 | 7.48 | 5.76 | 1.16 |
| **3** | 2249.09 | 2318.94 | 2278.94 | 19.31 | **-0.66** | 4.88 | 7.78 | 5.55 | 0.92 |
| **4** | 2561.86 | 2710.75 | 2605.43 | 43.26 | **1.30** | 4.27 | 5.53 | 4.89 | 0.37 |
| **5** | 2259.14 | 2320.20 | 2293.29 | 19.46 | **-0.87** | 4.67 | 7.92 | 5.49 | 0.95 |
| **6** | 2249.29 | 2291.62 | 2261.09 | 12.83 | **-0.56** | 3.72 | 5.63 | 5.07 | 0.57 |
| **7** | 2391.27 | 2500.45 | 2452.33 | 30.52 | **-0.86** | 7.50 | 8.32 | 7.82 | 0.25 |
| **8** | 2133.23 | 2223.38 | 2169.04 | 27.86 | **-4.77** | 7.63 | 10.30 | 8.73 | 0.83 |
| **9** | 2116.56 | 2173.60 | 2139.16 | 19.13 | **-1.51** | 8.62 | 13.20 | 9.38 | 1.43 |
| **10** | 2406.86 | 2456.27 | 2436.43 | 16.71 | **0.83** | 6.93 | 10.72 | 8.19 | 1.09 |
| **11** | 2151.52 | 2181.77 | 2169.24 | 12.15 | **-4.93** | 6.30 | 8.45 | 7.87 | 0.60 |
| **12** | 2100.29 | 2124.82 | 2112.60 | 9.62 | **-1.76** | 6.43 | 8.82 | 8.22 | 0.86 |
| **13** | 2386.49 | 2491.41 | 2423.94 | 27.93 | **-5.22** | 9.27 | 11.25 | 10.57 | 0.57 |
| **14** | 2223.47 | 2304.39 | 2258.76 | 22.69 | **-5.59** | 8.55 | 16.08 | 12.35 | 1.86 |
| **15** | 2208.77 | 2278.24 | 2232.98 | 20.27 | **-2.27** | 11.22 | 13.58 | 12.84 | 0.90 |
| **16** | 2374.48 | 2486.13 | 2426.60 | 32.73 | **-7.03** | 7.57 | 13.07 | 10.92 | 1.40 |
| **17** | 2209.06 | 2333.15 | 2266.24 | 37.47 | **-4.49** | 10.32 | 18.22 | 12.37 | 2.12 |
| **18** | 2197.09 | 2272.92 | 2230.88 | 28.06 | **-2.52** | 12.72 | 20.37 | 13.91 | 2.29 |

**Table C2.** Summarised numerical results of QL (E-greedy exploration)

| | cmax | | | | | Cpu (sec) | | | |
|---|---|---|---|---|---|---|---|---|---|
| **Inst.** | min | avg | max | stdev | %RE | min | avg | max | stdev |
| **1** | 2761 | 2825.4 | 2859 | 30.79 | **7.93** | 465.99 | 473.95 | 477.40 | 3.47 |
| **2** | 2348 | 2405.3 | 2441 | 27.98 | **3.64** | 519.67 | 531.84 | 550.47 | 9.88 |
| **3** | 2303 | 2318.6 | 2335 | 11.12 | **2.40** | 567.51 | 581.34 | 614.33 | 13.36 |
| **4** | 2775 | 2824.5 | 2861 | 23.21 | **8.32** | 466.72 | 472.68 | 481.43 | 4.34 |
| **5** | 2377 | 2392.3 | 2412 | 11.12 | **5.22** | 519.37 | 528.31 | 535.27 | 5.96 |
| **6** | 2291 | 2308.3 | 2320 | 10.20 | **1.85** | 571.53 | 579.57 | 596.36 | 7.73 |
| **7** | 2634 | 2664.8 | 2692 | 15.79 | **10.15** | 1133.34 | 1144.73 | 1156.16 | 7.08 |
| **8** | 2221 | 2268.2 | 2294 | 22.91 | **4.11** | 1308.40 | 1322.55 | 1335.62 | 9.32 |
| **9** | 2117 | 2144.6 | 2172 | 15.57 | **0.02** | 1541.86 | 1569.75 | 1641.03 | 29.93 |
| **10** | 2651 | 2680.6 | 2700 | 17.48 | **10.14** | 1140.30 | 1165.73 | 1228.40 | 29.80 |
| **11** | 2244 | 2268.6 | 2300 | 21.94 | **4.30** | 1304.42 | 1326.15 | 1343.65 | 12.29 |
| **12** | 2113 | 2138.6 | 2163 | 15.84 | **0.61** | 1544.38 | 1580.80 | 1607.12 | 21.29 |
| **13** | 2633 | 2682.4 | 2714 | 20.49 | **10.33** | 2072.75 | 2112.61 | 2177.06 | 40.24 |
| **14** | 2294 | 2316.7 | 2337 | 12.86 | **3.17** | 2475.40 | 2547.74 | 2620.43 | 42.83 |
| **15** | 2231 | 2242.7 | 2250 | 8.03 | **1.01** | 3103.58 | 3140.50 | 3197.72 | 27.77 |
| **16** | 2620 | 2663.1 | 2705 | 25.20 | **10.34** | 2079.09 | 2119.27 | 2243.16 | 50.43 |
| **17** | 2277 | 2311.6 | 2340 | 16.93 | **3.08** | 2488.37 | 2519.28 | 2566.77 | 25.39 |
| **18** | 2215 | 2237.6 | 2246 | 10.09 | **0.82** | 3073.96 | 3125.34 | 3156.40 | 28.77 |

**Table C3.** Summarised numerical results of QL (Pareto based selection)

| | cmax | | | | | Cpu (sec) | | | |
|---|---|---|---|---|---|---|---|---|---|
| **Inst.** | min | avg | max | stdev | %RE | min | avg | max | stdev |
| **1** | 2690 | 2737 | 2758 | 20.88 | **5.16** | 612.18 | 617.54 | 635.98 | 7.13 |
| **2** | 2363 | 2378.2 | 2393 | 10.60 | **4.30** | 696.74 | 706.45 | 717.51 | 7.02 |
| **3** | 2267 | 2278.3 | 2288 | 7.24 | **0.80** | 776.18 | 790.21 | 821.80 | 13.07 |
| **4** | 2690 | 2728.9 | 2758 | 21.02 | **5.00** | 607.24 | 614.94 | 624.78 | 4.58 |
| **5** | 2340 | 2358.8 | 2373 | 12.35 | **3.58** | 693.23 | 705.00 | 718.05 | 8.36 |
| **6** | 2265 | 2273.5 | 2283 | 5.19 | **0.70** | 784.32 | 795.39 | 820.38 | 11.29 |
| **7** | 2608 | 2630.7 | 2654 | 16.79 | **9.06** | 1467.47 | 1489.77 | 1545.28 | 23.12 |
| **8** | 2193 | 2221.7 | 2238 | 13.57 | **2.80** | 1865.85 | 1893.54 | 1961.12 | 31.83 |
| **9** | 2118 | 2125.6 | 2139 | 6.43 | **0.07** | 3064.72 | 3113.94 | 3164.81 | 32.36 |
| **10** | 2621 | 2634.5 | 2654 | 10.65 | **8.90** | 1467.88 | 1498.88 | 1604.24 | 40.20 |
| **11** | 2212 | 2231.1 | 2249 | 14.46 | **2.81** | 1853.69 | 1889.91 | 1912.45 | 18.72 |
| **12** | 2103 | 2115 | 2125 | 6.94 | **0.13** | 3070.54 | 3112.52 | 3206.32 | 39.16 |
| **13** | 2566 | 2629.2 | 2661 | 28.09 | **7.52** | 2853.08 | 2916.46 | 2966.89 | 28.67 |
| **14** | 2284 | 2300.7 | 2312 | 9.26 | **2.72** | 4830.86 | 4919.46 | 5089.88 | 69.86 |
| **15** | 2216 | 2226 | 2237 | 6.34 | **0.33** | 9540.91 | 9637.62 | 9728.62 | 59.62 |
| **16** | 2611 | 2636.3 | 2652 | 11.58 | **9.96** | 2890.07 | 2946.05 | 2977.20 | 32.36 |
| **17** | 2273 | 2288.3 | 2301 | 8.59 | **2.89** | 4882.65 | 4943.25 | 5016.07 | 39.33 |
| **18** | 2216 | 2220.7 | 2228 | 4.40 | **0.86** | 9571.43 | 9668.83 | 9804.67 | 62.30 |

**Table C4.** Number of states recorded and number of full solutions (E-greedy)

| | #States | | | | #Solutions | | | |
|---|---|---|---|---|---|---|---|---|
| **Inst.** | min | avg | max | stdev | min | avg | max | stdev |
| **1** | 7,081,810 | 7,088,903 | 7,095,802 | 4522.11 | 16580 | 16759.6 | 16872 | 98.96 |
| **2** | 7,283,528 | 7,288,098 | 7,291,856 | 3071.16 | 16594 | 16658.1 | 16729 | 48.93 |
| **3** | 7,440,054 | 7,444,102 | 7,453,290 | 5041.33 | 16565 | 16647 | 16699 | 39.51 |
| **4** | 7,092,481 | 7,097,748 | 7,099,952 | 2462.08 | 16239 | 16315.4 | 16404 | 51.13 |
| **5** | 7,297,152 | 7,300,983 | 7,304,096 | 2341.13 | 16567 | 16637.7 | 16735 | 62.85 |
| **6** | 7,444,526 | 7,446,749 | 7,450,988 | 2044.10 | 16470 | 16555.4 | 16697 | 67.98 |
| **7** | 10,472,253 | 10,475,147 | 10,480,199 | 2201.26 | 17576 | 17760 | 17833 | 79.73 |
| **8** | 10,904,277 | 10,911,508 | 10,926,435 | 9647.09 | 16245 | 16649.5 | 16853 | 224.67 |
| **9** | 11,285,969 | 11,288,081 | 11,289,492 | 1320.03 | 15825 | 15894.5 | 15968 | 48.10 |
| **10** | 10,488,404 | 10,489,985 | 10,491,699 | 1176.08 | 17645 | 17739.9 | 17838 | 59.29 |
| **11** | 10,910,645 | 10,919,882 | 10,932,239 | 8052.80 | 16416 | 16762.2 | 16925 | 156.27 |
| **12** | 11,287,531 | 11,288,653 | 11,290,423 | 979.53 | 15929 | 15992.4 | 16065 | 43.70 |
| **13** | 13,933,743 | 13,936,060 | 13,939,145 | 1738.92 | 16842 | 16934.5 | 16987 | 45.83 |
| **14** | 14,623,989 | 14,625,353 | 14,628,781 | 1428.32 | 15207 | 15265.1 | 15328 | 43.54 |
| **15** | 15,055,257 | 15,061,309 | 15,065,003 | 4106.29 | 14867 | 14996.5 | 15067 | 66.80 |
| **16** | 13,947,845 | 13,949,420 | 13,950,985 | 916.84 | 16783 | 16886.9 | 16966 | 58.45 |
| **17** | 14,638,289 | 14,641,043 | 14,646,234 | 2079.71 | 15155 | 15226.3 | 15303 | 53.32 |
| **18** | 15,054,653 | 15,059,247 | 15,061,916 | 2346.17 | 14664 | 14744.8 | 14841 | 58.32 |

**Table C5.** Number of states recorded and number of full solutions (Pareto based selection)

| | #States | | | | #Solutions | | | |
|---|---|---|---|---|---|---|---|---|
| **Inst.** | min | avg | max | stdev | min | avg | max | stdev |
| **1** | 7,141,300 | 7,141,936 | 7,142,864 | 426.86 | 15899 | 15972.4 | 16064 | 46.94 |
| **2** | 7,325,502 | 7,326,657 | 7,327,519 | 534.29 | 15346 | 15423 | 15501 | 46.79 |
| **3** | 7,532,704 | 7,533,429 | 7,533,912 | 448.02 | 15166 | 15233 | 15315 | 41.09 |
| **4** | 7,138,007 | 7,138,949 | 7,140,121 | 689.38 | 15982 | 16037.9 | 16103 | 41.04 |
| **5** | 7,338,947 | 7,339,526 | 7,340,119 | 293.86 | 15268 | 15341.8 | 15403 | 41.13 |
| **6** | 7,536,987 | 7,537,799 | 7,539,122 | 701.97 | 15258 | 15315.3 | 15360 | 34.49 |
| **7** | 10,512,697 | 10,513,280 | 10,513,856 | 372.99 | 15698 | 15825.7 | 15917 | 61.66 |
| **8** | 10,950,555 | 10,951,940 | 10,953,200 | 741.25 | 15077 | 15164 | 15263 | 57.15 |
| **9** | 11,381,526 | 11,382,209 | 11,383,124 | 451.32 | 15042 | 15081.4 | 15131 | 30.40 |
| **10** | 10,521,421 | 10,522,035 | 10,522,769 | 469.24 | 15725 | 15798.8 | 15876 | 47.10 |
| **11** | 10,958,767 | 10,959,262 | 10,960,060 | 377.81 | 15074 | 15155.7 | 15274 | 60.05 |
| **12** | 11,385,844 | 11,386,452 | 11,386,946 | 388.40 | 14947 | 15045.9 | 15147 | 64.66 |
| **13** | 14,006,332 | 14,007,681 | 14,008,667 | 650.63 | 16547 | 16599.9 | 16644 | 32.91 |
| **14** | 14,766,383 | 14,767,438 | 14,767,918 | 512.24 | 15523 | 15553.9 | 15589 | 24.68 |
| **15** | 15,158,021 | 15,158,636 | 15,159,331 | 342.81 | 14909 | 15027.6 | 15133 | 59.51 |
| **16** | 14,023,963 | 14,025,248 | 14,026,727 | 883.05 | 16497 | 16552.3 | 16677 | 52.41 |
| **17** | 14,775,737 | 14,776,800 | 14,777,947 | 764.70 | 15410 | 15479.7 | 15558 | 47.71 |
| **18** | 15,160,788 | 15,161,488 | 15,161,934 | 362.28 | 14938 | 15025.4 | 15147 | 57.21 |